\documentclass[%
  a4paper,
  onecolumn,
  colorlinks,
]{preprint}

\usepackage[english]{babel}
\usepackage[T1]{fontenc}

\usepackage{amsmath}
\usepackage{amssymb}
\usepackage{amsthm}
\usepackage{mathtools}
\usepackage{stmaryrd}
\SetSymbolFont{stmry}{bold}{U}{stmry}{m}{n}
\usepackage{mathrsfs}

\usepackage{graphicx}
\usepackage{xcolor}
\usepackage{subcaption}
\usepackage{tikz}
\usepackage{pgfplots}
\usepackage{pgfplotstable}
\usepackage[linesnumbered, lined, algoruled]{algorithm2e}

\usepackage{todonotes}
\makeatletter
\DeclareOldFontCommand{\rm}{\normalfont\rmfamily}{\mathrm}
\makeatother

\pgfplotsset{compat = 1.18}
\usepgfplotslibrary{external}

\tikzset{external/system call = {%
    pdflatex \tikzexternalcheckshellescape
    -halt-on-error
    -interaction=batchmode
    -jobname "\image" "\texsource"}}
\tikzexternaldisable

\theoremstyle{plain}

\theoremstyle{definition}
\newtheorem{remark}{Remark}
\newtheorem{definition}{Definition}
\newtheorem{proposition}{Proposition}

\newcommand{\R}{\ensuremath{\mathbb{R}}}

\newcommand{\bA}{\ensuremath{\boldsymbol{A}}}
\newcommand{\bB}{\ensuremath{\boldsymbol{B}}}
\newcommand{\bC}{\ensuremath{\boldsymbol{C}}}
\newcommand{\bD}{\ensuremath{\boldsymbol{D}}}

\newcommand{\bG}{\ensuremath{\boldsymbol{G}}}
\newcommand{\bH}{\ensuremath{\boldsymbol{H}}}
\newcommand{\bI}{\ensuremath{\boldsymbol{I}}}

\newcommand{\bL}{\ensuremath{\boldsymbol{L}}}
\newcommand{\bM}{\ensuremath{\boldsymbol{M}}}

\newcommand{\bP}{\ensuremath{\boldsymbol{P}}}
\newcommand{\bQ}{\ensuremath{\boldsymbol{Q}}}
\newcommand{\bR}{\ensuremath{\boldsymbol{R}}}
\newcommand{\bS}{\ensuremath{\boldsymbol{S}}}

\newcommand{\bU}{\ensuremath{\boldsymbol{U}}}
\newcommand{\bV}{\ensuremath{\boldsymbol{V}}}
\newcommand{\bW}{\ensuremath{\boldsymbol{W}}}
\newcommand{\bX}{\ensuremath{\boldsymbol{X}}}
\newcommand{\bY}{\ensuremath{\boldsymbol{Y}}}
\newcommand{\bZ}{\ensuremath{\boldsymbol{Z}}}

\newcommand{\bSigma}{\ensuremath{\boldsymbol{\Sigma}}}
\newcommand{\bOmega}{\ensuremath{\boldsymbol{\Omega}}}

\newcommand{\ba}{\ensuremath{\boldsymbol{a}}}
\newcommand{\bb}{\ensuremath{\boldsymbol{b}}}

\newcommand{\be}{\ensuremath{\boldsymbol{e}}}

\newcommand{\bm}{\ensuremath{\boldsymbol{m}}}
\newcommand{\bp}{\ensuremath{\boldsymbol{p}}}
\newcommand{\br}{\ensuremath{\boldsymbol{r}}}
\newcommand{\bs}{\ensuremath{\boldsymbol{s}}}

\newcommand{\bu}{\ensuremath{\boldsymbol{u}}}

\newcommand{\bx}{\ensuremath{\boldsymbol{x}}}
\newcommand{\by}{\ensuremath{\boldsymbol{y}}}

\newcommand{\bbeta}{\ensuremath{\boldsymbol{\beta}}}
\newcommand{\bell}{\ensuremath{\boldsymbol{\ell}}}
\newcommand{\bxi}{\ensuremath{\boldsymbol{\xi}}}
\newcommand{\blambda}{\ensuremath{\boldsymbol{\lambda}}}

\newcommand{\bTA}{\ensuremath{\boldsymbol{\mathcal{A}}}}
\newcommand{\bTB}{\ensuremath{\boldsymbol{\mathcal{B}}}}
\newcommand{\bTC}{\ensuremath{\boldsymbol{\mathcal{C}}}}
\newcommand{\bTG}{\ensuremath{\boldsymbol{\mathcal{G}}}}
\newcommand{\bTH}{\ensuremath{\boldsymbol{\mathcal{H}}}}
\newcommand{\bTM}{\ensuremath{\boldsymbol{\mathcal{M}}}}
\newcommand{\bTR}{\ensuremath{\boldsymbol{\mathcal{R}}}}

\newcommand{\bTV}{\ensuremath{\boldsymbol{\mathcal{V}}}}
\newcommand{\bTW}{\ensuremath{\boldsymbol{\mathcal{W}}}}
\newcommand{\bTX}{\ensuremath{\boldsymbol{\mathcal{X}}}}
\newcommand{\bTY}{\ensuremath{\boldsymbol{\mathcal{Y}}}}
\newcommand{\bTZ}{\ensuremath{\boldsymbol{\mathcal{Z}}}}

\newcommand{\trans}{\ensuremath{\mkern-1.5mu\mathsf{T}}}
\newcommand{\krp}{\odot}

\newcommand{\matlab}{\mbox{MATLAB}}
\newcommand{\plotfontsize}{\small}

\definecolor{matlabblue}{HTML}{0072BD}
\definecolor{matlabgreen}{HTML}{77AC30}
\definecolor{matlaborange}{HTML}{D95319}

\tikzset{Kron_Style_Bar/.style={
    matlabblue, 
    solid, 
    line width=2pt,
    mark=none,            
    fill=matlabblue!100
}}

\tikzset{KRP_Style_Bar/.style={
    matlaborange, 
    solid, 
    line width=2pt, 
    mark=none,                
    mark options={solid},
    fill=matlaborange!100 
}}

\tikzset{Naive_Style_Bar/.style={
    matlabgreen, 
    solid, 
    line width=2pt, 
    mark=none,         
    mark options={solid},
    fill=matlabgreen!100
}}

\tikzset{Naive_Style_line/.style={
    matlabgreen,
    solid,
    line width=2pt,
}}
    
\tikzset{Kron_Style_line/.style={
    matlabblue,
    solid,
    line width=2pt,
}}

\tikzset{KRP_Style_line/.style={
    matlaborange,
    solid,
    line width=2pt,
}}

\begin{document}


\title{Efficient Sketching-Based Summation of Tucker Tensors}

\author[1]{\mbox{Rudi Smith}}
\author[1]{\mbox{Mirjeta Pasha}}
\author[2]{\mbox{Andrés Galindo-Olarte}}
\author[3]{\mbox{Hussam Al Daas}}
\author[4]{\mbox{Grey Ballard}}
\author[5]{\mbox{Joseph Nakao}}
\author[6]{\mbox{Jing-Mei Qiu}}
\author[7]{\mbox{William Taitano}}

\affil[1]{Department of Mathematics, Virginia Tech, Blacksburg, VA 24061, USA.\authorcr
       \email{smithrgh@vt.edu}, \orcid{0000-0002-6592-2482}\authorcr
       \email{mpasha@vt.edu}, \orcid{0000-0003-4249-2421}}

\affil[2]{Oden Institute for Computational Engineering and Sciences, The University of Texas at Austin, Austin, TX 78712, USA.\authorcr
       \email{afgalindo@utexas.edu}, \orcid{0009-0007-1551-0889}}

\affil[3]{Computational Mathematics Group, Rutherford Appleton Laboratory, Chilton, Oxfordshire OX11 0QX, UK.\authorcr
       \email{hussam.al-daas@stfc.ac.uk}, \orcid{0000-0001-9355-4042}}

\affil[4]{Department of Computer Science, Wake Forest University, Winston Salem, NC 27109, USA.\authorcr
       \email{ballard@wfu.edu}, \orcid{0000-0003-1557-8027}}

\affil[5]{Department of Mathematics and Statistics, Swarthmore College, Swarthmore, PA 19081, USA.\authorcr
       \email{jnakao1@swarthmore.edu}, \orcid{0009-0008-6589-4013}}
\affil[6]{Department of Mathematical Sciences, University of Delaware, Newark, DE 19716, USA.\authorcr
       \email{jingqiu@udel.edu}, \orcid{0000-0002-3462-188X}}
\affil[7]{Los Alamos National Laboratory, Los Alamos, NM 87545, USA.\authorcr
       \email{taitano@lanl.gov}, \orcid{0000-0002-2369-0935}}

\shorttitle{Efficient Sketching-Based Summation of Tucker Tensors}
\shortauthor{R. Smith, M. Pasha, A. Galindo-Olarte, H. Al Daas, G. Ballard, J. Nakao, J. Qiu, W. Taitano}
\shortdate{2026-03-13}
\shortinstitute{}



\abstract{
We present efficient, sketching-based methods for the summation of tensors in Tucker format. Leveraging the algebraic structure of Khatri-Rao and Kronecker products, our approach enables compressed arithmetic on Tucker tensors while controlling rank growth and computational cost. The proposed sketching framework avoids the explicit formation of large intermediate tensors, instead operating directly on the factor matrices and core tensors to produce accurate low-rank approximations of tensor sums. Furthermore, we analyze the computational complexity and the theoretical approximation properties of the proposed methodology. Numerical experiments demonstrate the effectiveness of our approach on four problems: two synthetic test cases, a parameter-dependent elliptic equation (commonly referred to as the \emph{cookie problem}) solved via GMRES, and a one-dimensional linear transport problem discretized via high-order discontinuous Galerkin methods, where repeated tensor summation arises as a core computational bottleneck. Across these examples, the sketching-based summation achieves substantial computational savings while preserving high accuracy relative to direct summation and re-compression.
}

\novelty{%
  This work introduces efficient sketching-based algorithms for the summation of tensors in the Tucker format, by exploiting the algebraic structure of Khatri-Rao and Kronecker products. Unlike existing tensor compression techniques that rely on costly full-format operations or repeated truncations, our approach performs summation directly within the compressed Tucker representation, thereby precluding the explicit formation of large intermediate tensors. The proposed framework combines randomized sketching with structured matrix products to achieve substantial computational savings while empirically maintaining approximation accuracy. We demonstrate the practical efficacy of the method on problems arising from parametric PDEs and high order discontinuous Galerkin discretizations of linear transport problems, where repeated tensor summation is a core computational bottleneck.
}

\maketitle

\section{Introduction}%
\label{sec:Introduction}

In large-scale scientific computing, tensor computations have become ubiquitous across a wide range of disciplines, including signal processing, computer vision, quantum computing, and the numerical solution of high-dimensional partial differential equations~\cite{TokSPetal25, PanKCetal21, BerLAetal25,  Kho14, TiaPKetal24}.
As the dimensionality of the data increases, the storage and computational cost associated with grid-based representations grow exponentially, a phenomenon known as the `curse of dimensionality'~\cite{Bel66}.
To address this challenge, low-rank tensor decompositions have emerged as essential tools, including the Canonical Polyadic (CP)~\cite{CarC70, Har70}, Tensor Train (TT)~\cite{Ose11}, and Tucker decomposition~\cite{Tuc66}. These formats enable the compression of massive datasets into compact factor representations by exploiting latent multilinear structure.
While low-rank formats effectively solve the storage problem for static data, performing arithmetic operations on compressed tensors introduces a new set of computational challenges.  A fundamental issue arises during the summation of tensors, an operation that is central to iterative algorithms such as Krylov subspace methods (e.g., \texttt{GMRES}) for solving linear systems~\cite{BucPR25}. 
Unlike scalar arithmetic, where storage is constant, the rank of a tensor sum is generally the sum of the ranks of its constituents, which consequently means that repeated addition leads to a rapid increase in rank, quickly negating the compression benefits of the decomposition. To maintain computational feasibility in tensorial settings, it is usually standard practice to employ a rounding (truncation) strategy; the tensors are explicitly summed, and the resulting high-rank object is immediately rounded to a target tolerance~\cite{AldBCetal23}.

For a Tucker decomposition, the size of the core tensor scales with the product of the ranks along each mode. If we sum $d$ tensors, each with uniform rank $r$, the formal sum yields a block-sparse core. However, standard truncation procedures require orthogonalizing the factor matrices, which destroys this sparsity and forces the intermediate core tensor to swell to a dense $(dr)^N$ array, where $N$ is the order of the tensor. In high-dimensional settings (e.g., $N \ge 3$) or applications involving many summands, this intermediate growth of the core tensor often exceeds available memory, even if the final truncated result is compact. One may decide to truncate more frequently to avoid this memory growth, but this comes at its own cost of performing expensive rounding operations, losing some accuracy, or even rank swelling issues. 

To address this challenge, we propose a randomized approach to summing Tucker tensors that avoids the explicit formation of high-rank intermediate structures when rounding. Building on core concepts from Randomized Numerical Linear Algebra (RandNLA), our method leverages the concept of sketching where one maps high-dimensional data to a lower-dimensional subspace. Various forms of sketching have already been successfully applied to static tensors~\cite{PhaP25, HasN23, AldBCetal23, MinLB24, WanTSetal15}, which provide a good theoretical underpinning to our summation problem. Notably,~\cite{KreP17} highlights a highly relevant approach for recompressing Hadamard products in the Tucker format. While the Hadamard product affects the core tensor differently than a sum, their central idea of exploiting structure using randomized algorithms with rank-1 (Khatri-Rao) matrices is fundamentally similar to our approach.
In this work, we introduce a novel summation algorithm that employs structure preserving sketching operators based on Kronecker and Khatri-Rao products. By applying these sketching operators directly to the factor matrices of the individual summands, we construct a low-dimensional surrogate of the sum's principal subspace without ever realizing the full, high-rank core tensor. This technique effectively merges the summation and truncation steps into a single procedure, maintaining a bounded rank and memory footprint throughout the computation. To ensure the accuracy of the approximation while avoiding the pitfalls of over-sketching, our approach introduces a novel explicit step to estimate effective ranks prior to sketching. For each mode $n$, we construct an energy-weighted concatenated factor block $\bV_{n}$ and compute its corresponding Gram matrix $\bB_{n} = \bV_{n}^{\trans} \bV_{n}$. The effective rank $\widetilde{r}_n$ is then estimated based on the spectral decay of $\bB_{n}$, truncating eigenvalues that fall below a relative threshold and adding an oversampling parameter $p_n$ to ensure robustness. Upon making this estimation, we employ the subrank heuristic from~\cite{MinLB24} specifically for the Kronecker case when reusing factors. While we assume factor reuse throughout this paper, we note that in certain applications—such as those involving highly skewed modes—it may be more appropriate to utilize a separate heuristic without reusing factors, as noted in~\cite{MinLB24}.
We demonstrate the efficacy and robustness of this approach through a series of numerical experiments. First, we validate the scaling properties of our algorithm on synthetic tensors, showing that our runtime scales favorably with the number of summands compared to deterministic methods, and we justify the frequency of applying the baseline rounding approach. Subsequently, we apply the method to the \emph{cookie problem}~\cite{Tob12}, a parametric elliptic PDE with geometric inclusions that serves as a good benchmark for high-dimensional solvers~\cite{AldBCetal23, BalG15,  AldBGetal25}. Finally, we demonstrate the practical utility of our method in the context of a high-dimensional linear transport problem solved via the nodal discontinuous Galerkin (NDG) method~\cite{HesW08}, a prototype problem that naturally arises in low-rank kinetic simulations.

The remainder of the manuscript is organized as follows.
In \Cref{sec:prelim}, we review the necessary tensor notation and establish the algebraic foundations of Tucker decompositions, including deterministic summation and sequential truncation algorithms.
These results are then used to develop the randomized sketching framework for Tucker tensor sums utilizing Kronecker and Khatri-Rao products in \Cref{sec:Randomized_Sketching_of_Tucker_Tensor_Sums}, where the estimation of appropriate sketch dimensions and the heuristic selection of effective subranks are also discussed. In \Cref{sec:Complexity} we analyze the theoretical asymptotic complexities of the baseline and randomized summation algorithms.
In \Cref{sec:Results}, the proposed methods are tested in numerical experiments verifying their effectiveness and computational efficiency in comparison to the standard deterministic approaches across synthetic benchmarks and high-dimensional physical simulations.
The paper is concluded in \Cref{sec:conclusions}.

\section{Tensor Preliminaries}
\label{sec:prelim}

We refer the reader to~\cite{BalK25} and~\cite{HalMT11, MarT20} for general tensor and RandNLA preliminaries underpinning the randomized approach to summing Tucker tensors proposed in this work. In \Cref{subsec:TensorIndexingReshapingOperations} we review fundamental concepts from~\cite{BalK25} for indexing and reshaping tensors, and for basic tensor operations. In \Cref{subsec:TuckerTensors} we present the requisite Tucker tensor preliminaries, which establish the baseline methodologies for numerical comparisons later in \Cref{sec:Results}.

\subsection{Tensor Notation and Algebraic Operations}
\label{subsec:TensorIndexingReshapingOperations}

Throughout this work, we adopt a notational convention in which scalars are denoted by lowercase letters ($x$), vectors by bold lowercase letters ($\bx$), matrices by bold uppercase letters ($\bX$), and tensors by bold calligraphic letters ($\bTX$). We consider an order-$N$ tensor $\bTX \in \R^{n_{1} \times n_{2} \times \dots \times n_{N}}$ accessible via $N$ indices $x_{i_{1} i_{2} \dots i_{N}}$, such that $1 \le i_{k} \le n_{k}$ for all modes $k=1, \dots, N$. Fundamental to the multilinear operations discussed herein are the structured matrix products used to construct and manipulate low-rank tensor representations. 

\begin{definition}[Kronecker and Khatri-Rao products]
For matrices $\bA \in \R^{m \times n}$ and $\bB \in \R^{p \times q}$, the Kronecker product $\bA \otimes \bB \in \R^{mp \times nq}$ is defined by 
\begin{equation}
(\bA \otimes \bB)_{(i-1)p+k, (j-1)q+\ell} = a_{ij} b_{k\ell}.
\end{equation}

When the constituent matrices share the same column dimension (e.g., $\bA \in \R^{m \times r}$ and $\bB \in \R^{p \times r}$), we define a column-wise variant of this operation, known as the Khatri-Rao product (KRP) $\bA \krp \bB \in \R^{mp \times r}$
\begin{align}
    \bA \krp \bB = \begin{bmatrix} \ba_{1} \otimes \bb_{1} & \ba_{2} \otimes \bb_{2} & \dots & \ba_{r} \otimes \bb_{r} \end{bmatrix}.
\end{align}
\end{definition}

To apply these linear algebraic concepts to high-dimensional data, we require reshaping operations that map tensors to matrix representations. We refer to the one-dimensional sections of a tensor $\bTX$, obtained by varying the $k$-th index while fixing all others, as \textit{mode-$k$ fibers}, denoted in \matlab{} notation by $\bTX(i_{1}, \dots, i_{k-1}, :, i_{k+1}, \dots, i_{N})$. By arranging these fibers as the columns of a matrix, we obtain the \textit{mode-$k$ unfolding} (or matricization), denoted by $\bX_{(k)} \in \R^{n_{k} \times N_{k}}$, where $N_{k} = \prod_{j \neq k} n_{j}$.  The mapping between the tensor element indices and the matrix column index $j$ typically follows a lexicographical ordering of the modes excluding $k$
\begin{align}
    j = 1 + \sum_{\substack{m=1 \\ m \neq k}}^{N} (i_{m} - 1) \prod_{\substack{l=1 \\ l \neq k}}^{m-1} n_{l}. \label{eq:unfolding_index}
\end{align}

This matricization allows us to formally define the contraction of a tensor with a matrix. The \textit{mode-$k$ tensor-times-matrix product} (TTM) of a tensor $\bTX$ with a matrix $\bA \in \R^{r_{k} \times n_{k}}$, denoted as $\bTY = \bTX \times_k \bA$, applies the linear transformation $\bA$ to each mode-$k$ fiber of the tensor.  Element-wise, the resulting tensor $\bTY \in \R^{n_{1} \times \dots \times n_{k-1}\times r_{k} \times n_{k+1}\times \dots \times n_{N}}$ is given by the summation
\begin{align}
    y_{i_{1} \dots j \dots i_{N}} = \sum_{i_{k}=1}^{n_{k}} x_{i_{1} \dots i_{k} \dots i_{N}} a_{j i_{k}}.
\end{align}
Computationally, it is highly advantageous to express this operation in terms of the unfolding, where the mode-$k$ product corresponds directly to the standard matrix multiplication of the linear transformation $\bA$ and the unfolding of the original tensor
\begin{align}
    \bY_{(k)} = \bA \bX_{(k)}.
\end{align}

\subsection{Tucker Tensors}
\label{subsec:TuckerTensors}

The Tucker decomposition approximates a high-dimensional tensor via a core tensor multiplied by factor matrices along each mode, effectively compressing the data by identifying the dominant subspaces for each dimension. This concept can be naturally extended from its three-mode origins~\cite{Tuc66} to general $N$-mode decompositions. Instead of operating on the high-dimensional dense tensor, this decomposition allows for the direct manipulation of the significantly smaller core tensor and factor matrices. This approach drastically reduces computational overhead and is readily implemented in standard software packages~\cite{KolB09, KreT14, KosPAetal19}.

Given an $N$-th order tensor $\bTX \in \R^{n_1 \times \dots \times n_N}$ and target ranks $(r_1, \dots, r_N)$, the Tucker decomposition approximates $\bTX$ with a core tensor $\bTG \in \R^{r_1 \times \dots \times r_N}$ and a set of factor matrices $\{\bU_k \in \R^{n_k \times r_k}\}_{k=1}^N$ via sequential tensor-times-matrix (TTM) operations
\begin{align}
 \bTX \approx \bTG \times_1 \bU_1 \times_2 \bU_2 \dots \times_N \bU_N. \label{Eq:TuckerDecompositionModeProd}
\end{align}
 
Element-wise, this relationship can be expressed as
\begin{align}
 x_{i_{1} \dots i_{N}} \approx \sum_{j_{1}=1}^{r_{1}} \dots \sum_{j_{N}=1}^{r_{N}} g_{j_{1} \dots j_{N}} u^{(1)}_{i_{1} j_{1}} \dots u^{(N)}_{i_{N} j_{N}}, \label{Eq:TuckerDecompositionElementwise}
\end{align}
where $u_{i_{k}j_{k}}^{(k)}$ is the $(i_k,j_k)$-th element of the matrix $\bU_k$. For brevity, we denote a Tucker decomposition of this form as $\bTX = \llbracket \bTG; \bU_{1}, \dots, \bU_{N} \rrbracket$. Many deterministic and randomized algorithms have been proposed for the formulation of this approximation~\cite{AndB98, DeLDV00a, KolB09, BucR24, MalB18, HasN23}; however, throughout the remainder of this paper, we operate under the assumption that such a decomposition is provided \textit{a priori}, rather than requiring explicit formulation. We now formalize the mathematical connections between the decomposition, its mode unfoldings, and structured matrix products. These identities are fundamental to the sketching algorithms developed in subsequent sections.

\begin{proposition}[Key identities for Tucker Tensors~\cite{BalG15}]
\label{prop:key_identities}
Let $\bTX \in \R^{n_1 \times \dots \times n_N}$ be represented by the exact Tucker decomposition $\bTX = \llbracket \bTG; \bU_{1}, \dots, \bU_{N} \rrbracket$, then the following identities hold
\begin{enumerate}
    \item \textbf{Tucker unfolding identity}: The mode-$k$ unfolding of $\bTX$ can be expressed as
    \begin{align}
        \bX_{(k)} = \bU_{k} \bG_{(k)} \left( \bU_{N} \otimes \dots \otimes \bU_{k+1} \otimes \bU_{k-1} \otimes \dots \otimes \bU_{1} \right)^\top. \text{} \label{eq:TuckerUnfoldingKronecker}
    \end{align}
    
    \item \textbf{Tucker MTTKRP}: For a set of matrices $\bA_j \in \R^{n_j \times s}$ ($j \neq k$), the Matricized Tensor Times Khatri-Rao Product (MTTKRP) applied to the mode-$k$ unfolding of $\bTX$ satisfies
    \begin{align}
        \bX_{(k)} \left( \bA_{N} \krp \dots \krp \bA_{k+1} \krp \bA_{k-1} \krp \dots \krp \bA_{1} \right) = \bU_{k} \bG_{(k)} \left( \bU_{N}^\top \bA_{N} \krp \dots \krp \bU_{1}^\top \bA_{1} \right).
    \end{align}
    This identity demonstrates that the MTTKRP can be computed highly efficiently by projecting the individual factor matrices, entirely avoiding the explicit formation of the full tensor or the full KRP. 
    
    \item \textbf{Tucker Multi-TTM}: For a set of matrices $\bB_j \in \R^{s_j \times n_j}$ ($j \neq k$), the multiple tensor-times-matrix product applied to $\bTX$ yields
    \begin{align}
        \bTX \times_{j \neq k} \bB_j = \bTG \times_k \bU_k \times_{j \neq k} (\bB_j \bU_j). \text{}
    \end{align}
    This serves as the natural operational kernel for Kronecker-structured projections acting directly on the Tucker format.
\end{enumerate}
\end{proposition}

We now establish the algebraic framework for computing linear combinations of tensors in the Tucker format, which serves as our baseline for summation spectral rounding. Let $\bTX, \bTY \in \R^{n_{1} \times \dots \times n_{N}}$ be two tensors given by their Tucker decompositions $\bTX = \llbracket \bTG_{\bTX}; \bU_{1}, \dots, \bU_{N} \rrbracket$ and $\bTY = \llbracket \bTG_{\bTY}; \bV_{1}, \dots, \bV_{N} \rrbracket$, with respective multilinear ranks $\{r^{\bTX}_{k}\}_{k=1}^{N}$ and $\{r^{\bTY}_{k}\}_{k=1}^{N}$. We define the operation \texttt{TuckerAxby} to compute the weighted sum $\bTZ = \alpha \bTX + \beta \bTY$. The factor matrices of the resulting tensor $\bTZ$ are formed by the column-wise concatenation of the corresponding factor matrices of $\bTX$ and $\bTY$,
\begin{align}
    \bW_{n} = \begin{bmatrix} \bU_{n} & \bV_{n} \end{bmatrix} \in \R^{n_{n}\times (r_{n}^{\bTX} + r_{n}^{\bTY})}.
\end{align}
The core tensor $\bTG_{\bTZ}$ is constructed as a super-diagonal block tensor, placing the scaled cores $\alpha \bTG_{\bTX}$ and $\beta \bTG_{\bTY}$ in disjoint blocks along the $N$-dimensional super-diagonal. This resulting core possesses dimensions $(r_{1}^{\bTX} + r_{1}^{\bTY}) \times \dots \times (r_{N}^{\bTX} + r_{N}^{\bTY})$, but is stored as a sparse tensor. The complete summation procedure is formalized in \Cref{alg:tucker_axby}.

\begin{algorithm}[t]
    \caption{\texttt{TuckerAxby}} \label{alg:tucker_axby}
    \KwIn{Scalars $\alpha, \beta$; Tucker tensors $\bTX, \bTY$}
    \KwOut{Tucker tensor $\bTZ = \alpha \bTX + \beta \bTY$}
    
    $\bTG_{\bTX} \gets \alpha \cdot \bTG_{\bTX}$\;
    $\bTG_{\bTY} \gets \beta \cdot \bTG_{\bTY}$\;
    
    $N \gets \text{ndims}(\bTG_{\bTX})$\;
    
    \For{$n \gets 1$ \KwTo $N$}{
        $\bW_{n} \gets [\bU_{n}, \bV_{n}]$ \tcp*{Concatenate factors}
    }
    
    $\bs^{\bTX} \gets \text{size}(\bTG_{\bTX})$\;
    $\bs^{\bTY} \gets \text{size}(\bTG_{\bTY})$\;
    $\bTG_{\bTZ} \gets \text{zeros}(\bs^{\bTX} + \bs^{\bTY})$\;
    
    Construct indices $\mathcal{I}_{\bTX}$ covering $1 \dots s^{\bTX}_{n}$ for each mode $n$\;
    Construct indices $\mathcal{I}_{\bTY}$ covering $s^{\bTX}_{n} + 1 \dots s^{\bTX}_{n} + s^{\bTY}_{n}$ for each mode $n$\;
    
    $\bTG_{\mathcal{Z}}(\mathcal{I}_{\bTX}) \gets \bTG_{\bTX}$\;
    $\bTG_{\mathcal{Z}}(\mathcal{I}_{\bTY}) \gets \bTG_{\bTY}$ \tcp*{Not explicitly formed; store in sparse format}
    
    \Return $\bTZ = \llbracket \bTG_{\bTZ}; \bW_{1}, \dots, \bW_{N} \rrbracket$\;
\end{algorithm} 

A consequence of the additive operation defined above is the growth of the multilinear rank; specifically, if $\bTX$ and $\bTY$ have ranks $\{r^{\bTX}_{j}\}_{j=1}^{N}, \{r^{\bTY}_{j}\}_{j=1}^{N}$, respectively, then the sum $\bTZ$ generally has ranks $\{\br_k^{\bTX} + \br_k^{\bTY}\}_{k=1}^N$. In iterative methods, repeated application of this addition causes an explosion in both storage and computational costs. Although storing the formal sum requires no extra memory, the densification of the core for rounding triggers an unsustainable spike in RAM usage. In particular, in applications where a large number of Tucker tensors are continually summed, such as kinetic equations~\cite{GalNPetal25}, the memory requirements for compression become prohibitive. Therefore, we require an implicit mechanism to compress tensor sums. 

As a baseline, we employ a notion of \texttt{TuckerRounding} (or truncation)~\cite{DonYQetal23}, where the tensor is approximated by spectral truncation of the core tensor. We note a small difference in our approach since we do not assume orthogonal factor matrices. The procedure, detailed in \Cref{alg:st_hosvd_rounding}, first orthogonalizes the factor matrices via \texttt{QR} decomposition to transfer non-orthogonality into the core. Subsequently, a Higher-Order Singular Value Decomposition (\texttt{HOSVD})~\cite{DeLDV00b} style truncation is applied to the updated core tensor utilizing relative thresholding to achieve reconstruction within a specified tolerance. In practice, we actually apply a Sequentially Truncated Higher-Order Singular Value Decomposition (\texttt{ST-HOSVD})~\cite{VanVM12} to the updated core tensor. Unlike the standard \texttt{HOSVD}, which handles full-scale intermediate tensors prior to the final truncation, the \texttt{ST-HOSVD} enforces truncation immediately following the processing of each mode. This controls the growth of intermediate core structures that can occur with high-dimensional data.

\begin{remark}\label{rmk:rounding}
    There is a delicate interplay between computational cost and memory usage when applying a rounding procedure. Rounding too early may yield superfluous or negligible compression, whereas rounding too late results in a significant increase in storage requirements due to intermediate rank growth. Determining the optimal juncture for this operation is thus a non-trivial challenge in application and should be heavily considered in the experimental setup. We discuss the specific choices made for our experiments in \Cref{sec:Results}.
\end{remark}

We also make use of other standard Tucker tensor operations in our code~\cite{supSmiPGetal26}, based on the TT code for the parameter-dependent elliptic \emph{cookie problem}~\cite{AldBCetal23} using the tensor toolbox~\cite{KolB09}.

\begin{algorithm}[t]
    \caption{\texttt{TuckerRounding} using \texttt{ST-HOSVD}} \label{alg:st_hosvd_rounding}
    \SetKwInOut{Input}{Input}
    \SetKwInOut{Output}{Output}
    \DontPrintSemicolon
    
    \Input{Tucker tensor $\bTX = \llbracket \bTG; \bU_{1}, \dots, \bU_{N} \rrbracket$, relative tolerance $\tau$}
    \Output{Truncated Tucker tensor $\widetilde{\bTX}$}

    \tcp{Phase 1: Orthogonalize factor matrices to concentrate energy in the core}
    \For{$n \gets 1$ \KwTo $N$}{
        $[\bQ, \bR] \gets \text{\texttt{QR}}(\bU_{n}, \text{`econ'})$\;
        $\bU_{n} \gets \bQ$\;
        $\bTG \gets \bTG \times_{n} \bR$\;
    }
    
    \tcp{Distribute allowable error loss equally across all $N$ modes}
    $\theta \gets \frac{\tau^2 \|\bTG\|_F^2}{N}$\;
    
    \tcp{Phase 2: Sequential Truncation (ST-HOSVD) \cite{VanVM12}}
    \For{$n \gets 1$ \KwTo $N$}{
        $\bG_{(n)} \gets \texttt{unfold}(\bTG, n)$\;
        $[\bL, \bSigma, \sim] \gets \texttt{SVD}(\bG_{(n)}, \text{`econ'})$\;
        
        $\bs \gets \text{diag}(\bSigma)$\;
        
        $r_n \gets \min \left\{ r \geq 1 \ \middle|\ \sum_{j=r+1}^{|\bs|} s_j^2 \leq \theta \right\}$\; 
        
        $\bP \gets \bL_{:, 1:r_n}$\;
        
        $\bU_{n} \gets \bU_{n} \bP$\;
        $\bTG \gets \bTG \times_n \bP^{\trans}$\;
    }

    \Return $\widetilde{\bTX} = \llbracket \bTG; \bU_{1}, \dots, \bU_{N} \rrbracket$\;
\end{algorithm}

\section{Randomized Sketching of Tucker Tensor Sums}
\label{sec:Randomized_Sketching_of_Tucker_Tensor_Sums}

One advantage of a randomized approach is that it can avoid the truncation dilemma discussed in \Cref{rmk:rounding} by performing rounding implicitly; effectively performing a summation and truncation step without the intermediate memory overhead. From \Cref{prop:key_identities}, we apply sketching methods to the larger dimensions $n_{1}, \dots, n_{N}$ of the factor matrices to effectively compress each factor matrix before reconstructing the resulting tensor using \texttt{ST-HOSVD}. In particular, we seek a sketching matrix that can exploit the Kronecker structure of \Cref{eq:TuckerUnfoldingKronecker} without explicitly forming the Kronecker product.

\subsection{Sketching Methods}

Sketching is a dimensionality reduction technique~\cite{Woo14} where the fundamental principle is to compress a high-dimensional low-rank input $\bA \in \R^{m \times n}$ ($m \gg n$) into a significantly smaller $\bOmega\bA \in \R^{r \times n}$ by applying a random linear map $\bOmega \in \R^{r \times m}$, chosen such that key geometric properties, such as Euclidean norms and subspace angles, are preserved with high probability. By substituting the original data with the compressed surrogate, expensive linear algebraic computations can be performed in the lower-dimensional space, yielding approximate solutions at a fraction of the computational cost.

In standard summation, the concatenation of factor matrices and cores leads to rapid rank accumulation. While sketching could theoretically compress these factors post summation, doing so would require realizing the intermediate large-scale matrices, thereby negating the memory benefits. Instead, we employ structure-exploiting sketching techniques using the KRP~\cite{CheW19, SunGLetal20} and Kronecker product~\cite{CheWY20, CheWY21}, which allow us to apply the sketching operator directly to the constituent factors before the products are ever explicitly formed, significantly reducing computational complexity and memory footprint. 

\subsubsection{Khatri-Rao Sketching}

The first approach uses sketching matrices based on the KRP. We construct the $k$-th sketching matrix $\bS_{k}^{\text{KR}} \in \R^{N_k \times s}$ as the KRP of individual standard Gaussian random matrices $\bOmega_{j} \in \R^{n_{j} \times s}$ for each mode $j \neq k$, and note that the target sketch size $s$ is the same for all $\bOmega_{j}$ so that the operation is well-defined
\begin{align}\label{eq:KRP_sketch}
 \bS_{k}^{\text{KRP}} = \mathbf{\Omega}_N \krp \dots \krp \mathbf{\Omega}_{k+1} \krp \mathbf{\Omega}_{k-1} \krp \dots \krp \mathbf{\Omega}_1 := \bigodot_{\substack{j=N\\j\neq k}}^1 \mathbf{\Omega}_j.
\end{align}
While it is theoretically possible to generate distinct random matrices for each factor in the $k$-th sketch using Equation \eqref{eq:KRP_sketch}, this approach imposes a substantial computational burden due to the volume of random number generation required. To mitigate this, we adopt a strategy from~\cite{MinLB24}, wherein we reuse the sketched components across iterations. One drawback of sketching using the KRP is that we might sometimes over-sketch one of the factors if $s \geq n_j$. While one might attempt to circumvent this by reducing $s$, such a reduction introduces the risk of insufficiently sketching the most critical factor matrix; we develop a selection criterion that works well in practice in \Cref{sec:Subrank_Selection}. KRP sketching advantageously allows us to choose $s$ more specifically and increase the overall sketch size in smaller increments than Kronecker sketching, but has a limitation that we have to sketch all modes with this size. 

Using the KRP, we present an efficient way to compute the sketched unfolding for a single Tucker tensor. Let $\bTX = \llbracket \bTG; \bU_{1}, \dots, \bU_{N} \rrbracket$ and $\bS_{k}^{\text{KRP}}$ be as in \Cref{eq:KRP_sketch}; then the sketched mode-$k$ unfolding $\bX_{(k)} \bS_{k}^{\text{KRP}}$ is given by
 \begin{align}
  \bX_{(k)} \bS_k^{\text{KRP}} = \bU_k \bG_{(k)} \left( \bigodot_{\substack{j=N\\j\neq k}}^1 (\bU_j^{\trans} \mathbf{\Omega}_j) \right), \label{eq:krp_sketch_result}
 \end{align}
by \Cref{prop:key_identities} and the key identity $(\bA \otimes \bB) (\bC \odot \bD) = (\bA\bC) \odot (\bB\bD)$. This is significant because it avoids large matrix computations and allows us to initially sketch the factor matrices individually. Instead of forming \mbox{$\left(\bU_{N} \otimes \dots \otimes \bU_{k+1} \otimes \bU_{k-1} \otimes \dots \otimes \bU_{1} \right)^{\trans}$} and $\bS_{k}^{\text{KRP}}$ explicitly, we compute smaller matrices $\bM_{j} = \bU_{j}^{\trans} \bOmega_{j} \in\R^{r_j \times s}$. Applying this to the sum of tensors $\bTC = \sum_{i=1}^d \bTX^{(i)}$, where $i$ indexes the $i$-th tensor, gives the formula for the sketched sum
\begin{align}
\label{eq:KRP_sketch_unfolding}
 \bY_{k} = \bC_{(k)} \bS_{k}^{\text{KR}} = \sum_{i=1}^d \left[ \bU_{k}^{(i)} \bG_{(k)}^{(i)} \left( \bigodot_{\substack{j=N\\j\neq k}}^1 ((\bU_{j}^{(i)})^{\trans} \bOmega_{j}) \right) \right].
\end{align}
This allows efficient computation of the sketched sum $\bY_k \in \R^{n_k \times s}$ where $s$ is chosen to be slightly larger than the expected rank of $\bC_{(k)}$ via some oversampling parameter $p$.

\subsubsection{Kronecker Sketching}

The second approach utilizes sketching matrices predicated on the Kronecker product. Let $\bOmega_{j} \in \R^{n_{j} \times s_{j}}$ be a random sketching matrix for mode $j$. In this framework, the sketch sizes $s_{j}$ can differ across modes; however, this implies that the column dimension of the overall sketching matrix $\bS_{k}^{\text{Kron}} \in \R^{N_k \times S_{k}}$ scales multiplicatively, with $S_k = \prod_{j \neq k} s_{j}$. The operator is defined as
\begin{align}
 \bS_{k}^{\text{Kron}} = \bOmega_{N} \otimes \dots \otimes \bOmega_{k+1} \otimes \bOmega_{k-1} \otimes \dots \otimes \bOmega_{1} := \bigotimes_{\substack{j=N\\j\neq k}}^1 \bOmega_{j}.
\end{align}
We proceed similarly to the KRP case to establish an efficient sketching procedure. Let $\bTX = \llbracket \bTG; \bU_{1}, \dots, \bU_{N} \rrbracket$ and $\bS_{k}^{\text{Kron}} = \bigotimes_{\substack{j=N\\j\neq k}}^1 \bOmega_{j}$, then the sketched mode-$k$ unfolding is given by
 \begin{align}
  \bX_{(k)} \bS_{k}^{\text{Kron}} = \bU_{k} \bG_{(k)} \left( \bigotimes_{\substack{j=N\\j\neq k}}^1 (\bU_{j}^\top \bOmega_{j}) \right), \label{eq:kron_sketch_result}
 \end{align}
by using \Cref{prop:key_identities} and the standard mixed-product property of Kronecker products $(\bA \otimes \bB) (\bC \otimes \bD) = (\bA\bC) \otimes (\bB\bD)$. This formulation guarantees efficient computation by restricting operations to the smaller projected matrices $\bM_{j} = \bU_{j}^\top \bOmega_{j} \in \R^{r_{j}\times s_{j}}$, which, when applied to the summation $\bTC = \sum_{i=1}^d \bTX^{(i)}$, yield the formula for the sketched sum
\begin{align}
\label{eq:Kron_sketch_unfolding}
 \bY_{k} = \bC_{(k)} \bS_{k}^{\text{Kron}} = \sum_{i=1}^d \left[ \bU_k^{(i)} \bG_{(k)}^{(i)} \left( \bigotimes_{\substack{j=N\\j\neq k}}^1 ((\bU_{j}^{(i)})^\top \bOmega_{j}) \right) \right].
\end{align}
Determining the optimal sketch sizes for both Kronecker and KRP sketching remains a non-trivial problem; consequently, we develop a heuristic subrank selection method based on the spectral properties of the concatenated factor matrices of the sum.

\subsection{Heuristic Subrank Selection}
\label{sec:Subrank_Selection}

A fundamental challenge that remains is the selection of optimal sketch dimensions that properly balance computational efficiency with approximation accuracy. To address this, we adapt the subrank selection heuristic proposed by Minster et al.~\cite{MinLB24}, incorporating a spectral analysis of the summand subspaces to prevent the over-estimation of the target rank. We briefly recall that our summation target rank is bounded by the sum of the ranks. To determine a tighter effective rank $\widetilde{\br} \in \mathbb{N}^N$, we analyze the spectral decay of the weighted factor matrices of the sum. For a sum of tensors $\bTC = \sum_{i=1}^d \bTX^{(i)}$, where $\bTX^{(i)} = \llbracket \bTG^{(i)}; \bU_{1}^{(i)}, \dots, \bU_{N}^{(i)} \rrbracket$, we construct the energy-weighted concatenation for mode $n$ as
\begin{align}
    \bV_{n} = \left[ \|\bTG^{(1)}\|_F \bU_n^{(1)}, \dots, \|\bTG^{(d)}\|_F \bU_n^{(d)} \right].
\end{align}
The effective rank $\widetilde{r}_n$ is then determined from the eigenvalues $\lambda_j$ of the Gram matrix $\bV_{n}^{\trans} \bV_{n}$ by finding the minimum rank required to preserve sufficient spectral energy. Specifically, $\widetilde{r}_n$ is chosen such that the sum of the discarded eigenvalues is bounded by a relative threshold of $\frac{\epsilon^2}{N} \sum_{j} \lambda_j$ for a user-defined tolerance $\epsilon$ and tensor order $N$.
We define the target rank for sketching as $\bell = \widetilde{\br} + \bp$, where $\bp$ is a small, constant, user-defined safety oversampling parameter. In the Kronecker case, to ensure the sketching operator preserves the structural information of the tensor, we must select sketch dimensions $\bs \in \mathbb{N}^{N}$ that satisfy the projection requirements for each mode unfolding. The Kronecker structure of the sketch requires that the product of the sketch dimensions complementary to mode $n$ must define a subspace large enough to embed the rank $\ell_n$, yielding the condition
\begin{align}
    \prod_{j \neq n} s_{j} \geq \ell_{n}, \quad \forall n \in \{1, \dots, N\}.
    \label{eq:sketch_condition}
\end{align}
To determine a minimal set of parameters satisfying this bound, we solve for the limiting case where equality holds, implying that the product $s_n \ell_{n}$ is constant across all modes. By enforcing this constraint, we derive the closed-form heuristic from~\cite{MinLB24} for the $i$-th sketch dimension
\begin{align}
    s_{i} = \left\lceil \frac{ \left( \prod_{j=1}^{N} \ell_{j} \right)^{\frac{1}{N-1}} }{ \ell_i } \right\rceil.
    \label{eq:heuristic_formula}
\end{align}
This formulation distributes the sketching capacity inversely proportionally to the mode rank $\ell_i$, ensuring that high-rank modes are compressed less aggressively than low-rank modes. The complete procedure, combining the effective rank estimation via the weighted Gram matrix with the heuristic selection, is summarized in \Cref{alg:effective_subrank}. 

Importantly, we note that for Tucker tensors exhibiting highly skewed mode dimensions, this specific selection criterion is not appropriate in the Kronecker case, a limitation also acknowledged by the original authors. However, KRP sketching does not have this limitation in the factor reuse case where we simply select the uniform sketching dimension as $s = \max_i \ell_i$, ensuring that, at a minimum, the maximum rank is adequately sketched. As previously noted, this approach may result in the oversketching of certain dimensions; nevertheless, we demonstrate that this strategy remains highly effective in practice. 

The concatenation operation on the factor matrices required here can be computationally expensive, as it scales with the sum of the ranks, which is inherently undesirable. In our experiments, we nonetheless utilize this formulation to establish that, under this selection criterion, the accuracy of the baseline method is strictly preserved. We do not elaborate on further optimizations of this step, as empirical observations already indicate substantial computational speedups. However, in practical applications, one may wish to employ approximation techniques for this selection process.

\begin{algorithm}[t]
    \caption{\texttt{EffectiveSubrankHeuristic}}
    \label{alg:effective_subrank}
    \SetKwInOut{Input}{Input}
    \SetKwInOut{Output}{Output}
    \DontPrintSemicolon
    
    \Input{Set of tensors $\{\bTX^{(i)}\}_{i=1}^d$ with core tensors $\{\bTG^{(i)}\}$ and factor matrices $\{\bU_n^{(i)}\}$, user-defined tolerance $\epsilon$, oversampling vector $\bp$}
    \Output{Sketch dimensions $\bs = (s_1, \dots, s_N)$}
    
    \tcp{Phase 1: Estimate Effective Ranks via Weighted Gram Matrix (\Cref{sec:Subrank_Selection})}
    \For{$n \gets 1$ \KwTo $N$}{

        $\bV_{n} \gets \left[ \|\bTG^{(1)}\|_F \bU_n^{(1)}, \dots, \|\bTG^{(d)}\|_F \bU_n^{(d)} \right]$\;
        
        $\bB_{n} \gets \bV_{n}^{\trans} \bV_{n}$\;
        
        $\blambda \gets \text{sort}(\max(0, \text{eig}(\bB_{n})), \text{'descend'})$
        $\tau \gets \frac{\epsilon^2}{N} \sum_{j} \lambda_j$ 
        $ \widetilde{r}_n \gets \min \left\{ r \geq 1 \ \middle|\ \sum_{j=r+1}^{|\blambda|} \lambda_j \leq \tau \right\}$\;
        
        $\ell_n \gets \widetilde{r}_n + p_n$ 
    }
    
    \tcp{Phase 2: Compute Heuristic Sketch Dimensions via \Cref{eq:heuristic_formula}}
    $P \gets \prod_{j=1}^N \ell_j$\;
    $C \gets P^{\frac{1}{N-1}}$\;
    
    \For{$i \gets 1$ \KwTo $N$}{
        $s_i \gets \left\lceil \frac{C}{\ell_i} \right\rceil$\;
    }
    
    \Return $\bs$\;
\end{algorithm}

\subsection{Sketching-Based Summation of Tucker Tensors}

We saw in \Cref{eq:KRP_sketch_unfolding,eq:Kron_sketch_unfolding} how to compute the sketched unfolding matrix $\bY_n = \bC_{(n)} \bS_n$, where $\bC_{(n)}$ is the mode-$n$ unfolding of the sum and $\bS_n$ is a sketching operator composed of the Kronecker product or KRP of smaller random matrices. Now to formally build our summation estimate we must employ an approach analogous to the Randomized \texttt{ST-HOSVD}~\cite{HalMT11, MinLB24}, but strictly adapted for the implicit summation structure. The process consists of four distinct stages, which we describe as follows: effective rank estimation, randomized range finding, core projection, and sequential truncation.

Firstly, prior to generating the random embeddings, it is critical to determine the target sketch dimension. We employ the effective rank heuristic detailed in \Cref{alg:effective_subrank}, which analyzes the singular value decay of the energy-weighted Gram matrix of the concatenated factors. This provides an estimate of the subspace dimension $\bell$ required to capture the significant energy of the sum. Utilizing this estimate, we calculate the subrank sketch dimensions $\bs$ via the heuristic in \Cref{eq:heuristic_formula} (Kronecker) or by $s = \max_i \ell_i$ (KRP). Secondly, we compute an orthonormal basis $\widehat{\bU}_n$ for each mode $n$ of the sum by orthogonalizing the sketch $\bY_n$ using the standard \texttt{QR} decomposition. Thirdly, utilizing the approximate factors $\{\widehat{\bU}_n\}_{n=1}^N$, the core tensor of the sum, $\widetilde{\bTG}$, is obtained by projecting the original sum onto this newly derived basis.  By linearity, this is equivalent to projecting each summand individually and summing the results
\begin{align}
    \widetilde{\bTG} = \sum_{i=1}^d \bTG^{(i)} \times_1 (\widehat{\bU}_1^{\trans} \bU_1^{(i)}) \times_2 \dots \times_N (\widehat{\bU}_N^{\trans} \bU_N^{(i)}).
\end{align}
This operation is highly computationally efficient, as it involves solely small, dense matrix multiplications and core summations. Finally, the resulting core $\widetilde{\bTG}$ typically assumes dimensions corresponding to the sketch size $\bs$, which inherently includes a safety oversampling buffer. To produce a compact output satisfying the user-specified tolerance $\epsilon$, we perform a final compression step using the \texttt{ST-HOSVD}. The factor matrices are subsequently updated, yielding the final approximate decomposition $\llbracket \widetilde{\bTG}; \widetilde{\bU}_1, \dots, \widetilde{\bU}_N \rrbracket$. The complete procedure is outlined in \Cref{alg:kron_summation,alg:krp_summation} for Kronecker and KRP sketching, respectively, but full implementations can be found in~\cite{supSmiPGetal26}.

\begin{algorithm}[p]
\caption{Kronecker Sketching-Based Summation of Tucker Tensors (\texttt{KronSum-Tucker})}
\label{alg:kron_summation}
\DontPrintSemicolon
\KwIn{Set of Tucker tensors $\{\bTX^{(i)}\}_{i=1}^d$ with factors $\{\bU_n^{(i)}\}_{n=1}^N$, user-defined tolerance $\epsilon$, oversampling vector $\bp$}
\KwOut{Approximate Tucker decomposition $\widetilde{\bTC} = \llbracket \widetilde{\bTG}; \widetilde{\bU}_1, \dots, \widetilde{\bU}_N \rrbracket$}
\BlankLine

\tcc{Stage 1: Effective Rank Estimation (\Cref{alg:effective_subrank})}
$\bs \gets \texttt{EffectiveSubrankHeuristic}(\{\bTX^{(i)}\}, \bp, \epsilon)$\;

\tcc{Stage 2: Randomized Range Finder via \Cref{eq:Kron_sketch_unfolding}}
Generate Gaussian test matrices $\{\bOmega_{n} \in \mathbb{R}^{n_n \times s_n}\}_{n=1}^N$\;
Initialize storage for sketched terms $\bZ_n^{(i)}$\;
\For{$i \gets 1$ \KwTo $d$}{
    \For{$n \gets 1$ \KwTo $N$}{
        $\bTV \gets \bTG^{(i)}$\;
        \For{$j \in \{1,\dots,N\} \setminus \{n\}$}{
            $\bM_j \gets (\bOmega_{j})^{\trans} \bU_j^{(i)}$\;
            $\bTV \gets \bTV \times_j \bM_j$\;
        }
        $\bZ_n^{(i)} \gets \bU_n^{(i)} \bV_{(n)}$\;
    }
}
\For{$n \gets 1$ \KwTo $N$}{
    Initialize accumulator $\bY_n \gets \mathbf{0}$\;
    \For{$i \gets 1$ \KwTo $d$}{
        $\bY_n \gets \bY_n + \bZ_n^{(i)}$\;
    }
    $[\widehat{\bU}_n, \sim] \gets \text{qr}(\bY_n, \text{`econ'})$\;
}

\tcc{Stage 3: Core Projection and Assembly}
Initialize intermediate core $\bTH \gets \mathbf{0}$\;
\For{$i \gets 1$ \KwTo $d$}{
    $\bTV \gets \bTX^{(i)} \times_1 \widehat{\bU}_1^{\trans} \times_2 \dots \times_N \widehat{\bU}_N^{\trans}$\;
    $\bTH \gets \bTH + \bTV$\;
}

\tcc{Stage 4: Final Compression (ST-HOSVD)}
$[\widetilde{\bTG}, \{\bV_n\}_{n=1}^N] \gets \texttt{ST-HOSVD}(\bTH, \epsilon)$\;

\For{$n \gets 1$ \KwTo $N$}{
    $\widetilde{\bU}_n \gets \widehat{\bU}_n \bV_n$\;
}

\Return $\llbracket \widetilde{\bTG}; \bU_1, \dots, \bU_N \rrbracket$\;
\end{algorithm}

\begin{algorithm}[p]
\caption{KRP Sketching-Based Summation of Tucker Tensors (\texttt{KRPSum-Tucker})}
\label{alg:krp_summation}
\DontPrintSemicolon
\KwIn{Set of Tucker tensors $\{\bTX^{(i)}\}_{i=1}^d$ with factors $\{\bU_n^{(i)}\}$, user-defined tolerance $\epsilon$, oversampling parameter $p$}
\KwOut{Approximate Tucker decomposition $\widetilde{\bTC} = \llbracket \widetilde{\bTG}; \widetilde{\bU}_1, \dots, \widetilde{\bU}_N \rrbracket$}
\BlankLine

\tcc{Stage 1: Effective Rank Estimation (\Cref{alg:effective_subrank})}
$s \gets \texttt{EffectiveSubrankHeuristic}(\{\bTX^{(i)}\}, p, \epsilon)$\;

\tcc{Stage 2: Randomized Range Finder via \Cref{eq:KRP_sketch_unfolding}}
Generate Gaussian test matrices $\{\bOmega_{n} \in \mathbb{R}^{n_n \times s}\}_{n=1}^N$\;
Initialize storage for sketched terms $\bZ_n^{(i)}$\;
\For{$i \gets 1$ \KwTo $d$}{
    \For{$n \gets 1$ \KwTo $N$}{
        \For{$j \in \{1,\dots,N\} \setminus \{n\}$}{
            $\bM_j \gets (\bU_j^{(i)})^{\trans} \bOmega_j$\;
        }
        $\bV_{(n)} \gets \bG_{(n)}^{(i)} \left( \bigodot_{\substack{j=N\\j\neq n}}^1 \bM_j \right)$\;
        $\bZ_n^{(i)} \gets \bU_n^{(i)} \bV_{(n)}$\;
    }
}
\For{$n \gets 1$ \KwTo $N$}{
    Initialize accumulator $\bY_n \gets \mathbf{0}$\;
    \For{$i \gets 1$ \KwTo $d$}{
        $\bY_n \gets \bY_n + \bZ_n^{(i)}$\;
    }
    $[\widehat{\bU}_n, \sim] \gets \text{qr}(\bY_n, \text{`econ'})$\;
}

\tcc{Stage 3: Core Projection and Assembly}
Initialize intermediate core $\bTH \gets \mathbf{0}$\;
\For{$i \gets 1$ \KwTo $d$}{
    $\bTV \gets \bTX^{(i)} \times_1 \widehat{\bU}_1^{\trans} \times_2 \dots \times_N \widehat{\bU}_N^{\trans}$\;
    $\bTH \gets \bTH + \bTV$\;
}

\tcc{Stage 4: Final Compression (ST-HOSVD)}
$[\widetilde{\bTG}, \{\bV_n\}_{n=1}^N] \gets \texttt{ST-HOSVD}(\bTH, \epsilon)$\;

\For{$n \gets 1$ \KwTo $N$}{
    $\bU_n \gets \widehat{\bU}_n \bV_n$\;
}

\Return $\llbracket \widetilde{\bTG}; \widetilde{\bU}_1, \dots, \widetilde{\bU}_N \rrbracket$\;
\end{algorithm}

\section{Computational Complexity}
\label{sec:Complexity}

To evaluate the computational advantages of the proposed sketching framework, we analyze the theoretical asymptotic complexities of the baseline and randomized summation algorithms. Let $N$ denote the order of the tensors, $I = \max_n n_n$ be the maximum mode dimension across all spatial or parameter domains, and $d$ be the number of summands. To maintain clarity in the asymptotic analysis, we assume each summand $\bTX^{(i)}$ possesses a uniform multilinear rank bounded by $r$. Furthermore, to achieve a target sketched subspace capacity of $r+p$, we define two distinct sketch dimension variables. For the Kronecker sketching method, we let $s$ denote the uniform \textit{mode-wise} sketch dimension bound, yielding a global sketch capacity of $s^{N-1} \approx r+p$. For the KRP sketching method, we let $m$ denote the \textit{global} sketch dimension, where $m \approx r+p$.

\texttt{TuckerAxby}: The deterministic approach operates in two distinct phases: forming the exact naive sum and subsequently compressing it. As detailed in \Cref{alg:tucker_axby}, the explicit summation requires concatenating the factor matrices for each of the $N$ modes. Appending $d$ factor matrices of dimension $I \times r$ creates a new factor matrix of size $I \times (dr)$, requiring $\mathcal{O}(N I d r)$ allocation operations. Simultaneously, the algorithm constructs an expanded core tensor by placing the individual $r \times \dots \times r$ cores along the super-diagonal of a new $N$-dimensional structure. Because this formal sum maintains a strictly block-sparse structure, instantiating it requires only $\mathcal{O}(d r^N)$ memory operations. 
The total cost is 
\begin{equation}
    C_{\rm TAxBy} = \mathcal{O}\left(N I d r + d r^N\right).
\end{equation}

\texttt{TuckerRounding}: Here we analyze the total cost of \Cref{alg:st_hosvd_rounding} that we denote by $C_{\rm TR}$.
Following \texttt{TuckerAxby}, \Cref{alg:st_hosvd_rounding} applies the Sequentially Truncated Higher-Order Singular Value Decomposition (\texttt{ST-HOSVD}) to the intermediate tensor. The first step involves orthogonalizing the concatenated factor matrices via an economy \texttt{QR} decomposition. Computing the \texttt{QR} decomposition of an $I \times (dr)$ matrix costs $\mathcal{O}(I (dr)^2)$ operations; repeating this across all $N$ modes yields $\mathcal{O}(N I d^2 r^2)$. 

The most severe computational bottleneck occurs next, when the algorithm absorbs the upper-triangular $\bR$ factors into the core tensor via tensor-times-matrix (TTM) products. It is this specific operation that destroys the block-sparse structure of the formal sum, forcing the core to densify into a massive $(dr) \times \dots \times (dr)$ array. This densification demands $\mathcal{O}((dr)^N)$ memory operations to allocate, and computing the sequence of $N$ dense TTM products costs $\mathcal{O}(N (dr)^{N+1})$.

Following this densification, \texttt{ST-HOSVD} matricizes the core tensor along a given mode, forming a matrix of size $(dr) \times (dr)^{N-1}$, and computes its Singular Value Decomposition (SVD). The SVD of an $m \times p$ matrix (where $m \le p$) scales as $\mathcal{O}(m^2 p)$, meaning this initial decomposition costs $\mathcal{O}((dr)^2 \cdot (dr)^{N-1}) = \mathcal{O}((dr)^{N+1})$. Because the core is immediately truncated along this mode, the subsequent SVDs for the remaining modes operate on strictly smaller tensors and are thus of lower asymptotic order. 

Consequently, the overall computational complexity for the deterministic sum-and-round approach is strictly dominated by the initial core densification step, yielding a total complexity of $\mathcal{O}(N I d^2 r^2 + N d^{N+1} r^{N+1})$. This severe high-degree polynomial dependence on the number of summands $d$ (specifically, scaling as $d^{N+1}$ compared to the desired $d^2$ scaling) mathematically isolates the intermediate rank explosion that renders standard methodology intractable for large-scale iterative summation. 
The total cost is 
\begin{equation}
    C_{\rm TR} = \mathcal{O}\left(N I d^2 r^2 + N d^{N+1} r^{N+1}\right).
\end{equation}

Conversely, the randomized approaches detailed in \Cref{alg:kron_summation} and \Cref{alg:krp_summation} fundamentally bypass the explicit formation of the massive $dr$-rank intermediate core tensor. Both sketching methods share an identical initialization cost for the effective rank estimation procedure outlined in \Cref{alg:effective_subrank}. This requires constructing the energy-weighted concatenated factors of size $I \times (dr)$ and forming their Gram matrices, costing $\mathcal{O}(I (dr)^2)$ per mode, totaling $\mathcal{O}(N I d^2 r^2)$ operations.

\texttt{KronSum-Tucker}: \Cref{alg:kron_summation} first projects the individual summand factor matrices using the random test matrices, costing $\mathcal{O}(d N I r s)$. However, the Kronecker product structure implies that the column dimension of the global sketching operator scales multiplicatively as $s^{N-1}$. When the algorithm accumulates the implicitly sketched unfoldings to find the randomized range, it must perform a \texttt{QR} decomposition on a matrix of size $I \times s^{N-1}$. This exacts a computational cost of $\mathcal{O}(I (s^{N-1})^2) = \mathcal{O}(I s^{2N-2})$. Including the cost of projecting the original core onto this new basis and applying the final \texttt{ST-HOSVD} to the resulting $s \times \dots \times s$ dense core, the total complexity for Kronecker sketching evaluates to $\mathcal{O}(N I d^2 r^2 + d N I r s^{N-1} + N I s^{2N-2} + N s^{N+1})$.While this effectively decouples the polynomial scaling from the number of summands $d$, the $\mathcal{O}(I s^{2N-2})$ orthogonalization cost highlights a structural sensitivity of Kronecker sketching. To maintain efficiency, the sketch dimension must be carefully chosen such that $s \approx (r+p)^{1/(N-1)}$. However, if highly skewed mode ranks or integer constraints force the global column dimension $s^{N-1}$ to overshoot the target rank, the sketched matrix becomes excessively wide, leading to severe computational penalties.
The total cost is 
\begin{equation}
    C_{\rm KST} = \mathcal{O}\left(N I d^2 r^2 + d N I r s^{N-1} + N I s^{2N-2} + N s^{N+1}\right).
\end{equation}

\texttt{KRPSum-Tucker}: \Cref{alg:krp_summation} projects the $d$ factor matrices into the uniform sketch dimension $m$ requiring $\mathcal{O}(d N I r m)$ operations. The critical computational advantage arises during the implicit sketching of the core tensor via the MTTKRP. Because the KRP shares a uniform column dimension, contracting the $r^N$ core with the sketched factors costs only $\mathcal{O}(m r^N)$ per mode for each of the $d$ summands, yielding $\mathcal{O}(d N m r^N)$. Furthermore, accumulating the sketched terms results in an $I \times m$ matrix; the subsequent \texttt{QR} decomposition to extract the orthonormal basis costs $\mathcal{O}(N I m^2)$ operations. The final core projection and \texttt{ST-HOSVD} compression steps are strictly bounded by $\mathcal{O}(d N r^N m + N m^{N+1})$. 

Summing these contributions, the total asymptotic complexity for \texttt{KRPSum-Tucker} is bounded by

\begin{equation}
    C_{\rm KRST} = \mathcal{O}\left(N I d^2 r^2 + d N m r^N + N I m^2 + N m^{N+1}\right).
\end{equation}

This derivation shows why the KRP approach is the most efficient and robust methodology evaluated in our numerical experiments when reusing factors. It explicitly eliminates the $\mathcal{O}(N I d^2 r^2 + d^{N+1} r^{N+1})$ bottleneck of the deterministic method. Furthermore, it bypasses the potentially prohibitive orthogonalization penalty of the Kronecker method. While Kronecker sketching attempts to balance its mode-wise sketch size such that $s^{N-1} \approx m$, integer constraints and skewed ranks frequently force $s^{N-1} \gg m$, causing Kronecker to drastically over-sample the target subspace. KRP sketching avoids this geometric constraint entirely by directly utilizing $m$ columns, ensuring that both the computational runtime and the memory overhead remain strictly linear with respect to the number of summands $d$, facilitating highly scalable summation for high-dimensional tensors.

\section{Numerical Experiments}
\label{sec:Results}

The reported synthetic numerical experiments in \Cref{subsec:Synthetic} were performed on a MacBook Air with 16\,GB of RAM and an Apple M2 processor running macOS Ventura version 13.4 with MATLAB 25.1.0.2973910 (R2025a) Update 1. The reported large-scale numerical experiments for the high dimensional elliptic PDE in \Cref{subsec:cookie_problem} and the linear conservation law problem in \Cref{subsec:Linear_conservation_law} were performed with computational resources provided by Advanced Research Computing (ARC) at Virginia Tech (see acknowledgments) using a system with 128\,GB of RAM running Linux (5.14.0-611.13.1.el9\_7.x86\_64) with MATLAB 25.2.0.2998904 (R2025b). The source codes, data, and results of the numerical experiments reported in
this section are available at~\cite{supSmiPGetal26}.

\subsection{Synthetic Examples}
\label{subsec:Synthetic}
\subsubsection{Sketching a Summation that is Low-Rank}

We start by verifying that the proposed methodologies behave as theoretically expected before scaling them to larger, more complex problems. To evaluate the scalability of the randomized algorithms compared to the deterministic `sum and round' approach, we vary the number of summands $d$ while holding the tensor dimensions constant. In each \texttt{TuckerAxby} operation, an increasingly large intermediate core is sparsely stored prior to a final truncation via \texttt{TuckerRounding}. These experiments concern a fourth-order tensor ($N=4$) with a fixed mode dimension of $n_k=4000$ for all $k \in \{1, \dots, 4\}$. The number of summands $d$ is varied across the set $\{20, 40, 60, 80, 100\}$. To test the capacity of the algorithms to identify latent low-rank structures within nominally high-rank sums, we construct synthetic data embedded in a shared low-rank subspace, ensuring there is no true mathematical rank growth. Specifically, for each mode $k$, we first generate a fixed orthogonal basis $\bQ_k \in \R^{4000 \times r_{\text{true}}}$ via the economy-size \texttt{QR} decomposition of a standard Gaussian random matrix, fixing the true latent rank at $r_{\text{true}} = 10$. Each of the $d$ summands $\bTX^{(i)}$ is subsequently generated as a rank-1 Tucker tensor (core rank $1 \times 1 \times 1 \times 1$), where the corresponding factor matrix $\bU^{(i)}_{k}$ for mode $k$ is constructed as a random linear combination of the columns of the fixed basis:
\begin{align}
    \bU^{(i)}_{k} = \bQ_k \bm^{(i)}_k,
\end{align}
where $\bm^{(i)}_k \in \R^{r_{\text{true}} \times 1}$ is a standard Gaussian random mixing vector. This structural construction guarantees that the factor matrices of the sum lie entirely within the column space of $\bQ_k$. Consequently, the true mathematical rank of the final evaluated tensor remains exactly $r_{\text{true}} = 10$, despite the nominal rank of the naive sum growing linearly with $d$.

We compare \texttt{TuckerAxby}+\texttt{TuckerRounding}, where the truncation process requires forming a massive dense intermediate core of size up to $100^4$, against \texttt{KronSum-Tucker} and \texttt{KRPSum-Tucker}. The randomized methods employ an oversampling parameter of $p=2$, capturing the true rank of the solution precisely, and all rounding operations utilize a strictly enforced relative tolerance of $\epsilon = 10^{-6}$. To mitigate stochastic variance, all relative error metrics and runtime recordings are averaged over $S=3$ independent trials. The computational runtime analysis, alongside the resulting approximation errors and output ranks, is summarized in \Cref{fig:Synthetic_Example} where the standard deterministic approach suffers as the number of summands $d$ increases due to the densification of the massively swollen intermediate core. In contrast, both randomized sketching algorithms successfully circumvent this dense intermediate rank explosion, maintaining a dramatically lower runtime across all tested values $d$. Furthermore, the accompanying table demonstrates that this computational efficiency does not come at the cost of accuracy; both \texttt{KronSum-Tucker} and \texttt{KRPSum-Tucker} achieve approximation errors near machine precision relative to $\texttt{TuckerAxby+TuckerRounding}$. Ultimately, these results confirm that the randomized sketching techniques efficiently and capture the latent low-rank structure without the prohibitive memory and time costs inherent to densifying the naive sequential sum, justifying application to larger problems.

\begin{figure}[t]
    \centering
    \begin{subfigure}{0.99\textwidth}
        \centering
  \tikzexternalenable%
  \tikzsetnextfilename{Synthetic_Runtime}%
  \begin{tikzpicture}[font = \plotfontsize\normalfont]
  \pgfplotstableread{graphics/data/Synthetic_Times.dat}\tableSyntheticRuntime

  \begin{axis}[
    ybar,
    width              = \linewidth,
    height             = 5cm,
    enlarge x limits   = 0.15,
    bar width          = 20pt,
    ymin               = 0.01,
    ymode              = log,
    log origin         = infty,
    xtick              = data,
    xticklabels from table = {\tableSyntheticRuntime}{d},
    x tick label style = {/pgf/number format/1000 sep=},
    xlabel             = {Number of Summands ($d$)},
    ylabel             = {Total Runtime (s)},
    ylabel style       = {yshift = -.3em},
    ymajorgrids        = true,
    scaled x ticks     = false,
  ]
    
    \addplot[Kron_Style_Bar, mark=square*, mark options={black, xshift=-22pt, yshift=-8pt, thick}] 
    table[x expr=\coordindex, y=Time_Kron] {\tableSyntheticRuntime};

    \addplot[KRP_Style_Bar, mark=*, mark options={black, xshift=0pt, yshift=-8pt}] 
    table[x expr=\coordindex, y=Time_KRP] {\tableSyntheticRuntime};

    \addplot[Naive_Style_Bar, mark=triangle*, mark options={black, xshift=22pt, yshift=-8pt}] 
    table[x expr=\coordindex, y=Time_Det] {\tableSyntheticRuntime};

  \end{axis}
\end{tikzpicture}%
  \tikzexternaldisable%

        \label{fig:Synthetic_Runtime}
    \end{subfigure}
    \begin{subfigure}{0.45\textwidth}
        \centering
  \tikzexternalenable%
  \tikzsetnextfilename{Synthetic_Ranks_Error}%
  \begin{tikzpicture}[font = \plotfontsize\normalfont]
    \node[inner sep=0pt] {\pgfplotstabletypeset[
        col sep=space,
        every row/.style={after row=\hline},
        every last row/.style={after row=\hline},
        every head row/.style={
            before row={
                \hline
                \multicolumn{1}{|c|}{} & \multicolumn{2}{c|}{Relative Approximation Error} \\ 
                \hline
            },
            after row=\hline
        },
        columns={d_Summands, Error_Kron, Error_KRP},
        columns/d_Summands/.style={
            column name={$d$}, 
            int detect,
            column type=|c| 
        },            
        columns/Error_Kron/.style={
            column name={Kron}, 
            sci, sci zerofill, precision=2,
            column type=c|
        },
        columns/Error_KRP/.style={
            column name={KRP}, 
            sci, sci zerofill, precision=2,
            column type=c|
        }
    ]{graphics/data/Synthetic_Errors.dat}};
\end{tikzpicture}%
  \tikzexternaldisable%
 
        \label{fig:Synthetic_Ranks_ERROR_Accuracy}
    \end{subfigure}

          \vspace{.5\baselineskip}
  \tikzexternalenable%
  \tikzsetnextfilename{Legend}%
  \begin{tikzpicture}[font = \plotfontsize\normalfont]
  \begin{axis}[%
    hide axis, 
    scale only axis,
    width          = .8\linewidth,
    height         = .2\linewidth,
    xmin           = 0,
    xmax           = 1,
    ymin           = 0,
    ymax           = 1,
    legend columns = 3,    
    legend style   = {
      at     = {(0.5,0.5)}, 
      anchor = center,
      draw   = none,       
      /tikz/every even column/.append style = {column sep = 0.5cm}
    }
  ]
    
    
    \addlegendimage{Kron_Style_line, mark=square*, mark options={black}}
    \addlegendentry{\texttt{KronSum-Tucker}}

    \addlegendimage{KRP_Style_line,mark=*, mark options={black}}
    \addlegendentry{\texttt{KRPSum-Tucker}}

    \addlegendimage{Naive_Style_line,mark=triangle, mark options={black}}
    \addlegendentry{\texttt{TuckerAxby}+\texttt{TuckerRounding}}

  \end{axis}
\end{tikzpicture}%
  \tikzexternaldisable%

  \caption{Performance comparison for the synthetic low-rank summation example: 
The runtime of the \texttt{TuckerAxby}+\texttt{TuckerRounding} method grows rapidly with the number of summands, whereas the randomized sketching approaches scale highly efficiently. 
At $d=100$, the randomized KRP method is nearly two orders of magnitude faster than the deterministic baseline. Crucially, this computational speedup does not compromise accuracy; both the Kronecker and KRP sketching methods maintain errors relative to the deterministic baseline near machine precision ($\sim 4 \cdot 10^{-15}$) across all tested summation sizes.}
    \label{fig:Synthetic_Example}
  
\end{figure}

\subsubsection{Avoiding Truncation at Every Summation Step}

To explicitly demonstrate the algorithmic nuance and inherent danger of truncating sequentially immediately following every \texttt{TuckerAxby} operation (referred to herein as \textit{eager rounding}), we construct a specific synthetic example that highlights two severe numerical vulnerabilities: intermediate rank swelling and irreversible accuracy loss driven by dynamic range shifts. When summands contain large-magnitude, opposing components that eventually cancel analytically, sequential truncation forces the \Cref{alg:tucker_axby} to perform computationally heavy singular value decompositions on artificially inflated intermediate states. Simultaneously, this process introduces the risk of permanently deleting smaller, physically relevant signals. 

To evaluate this phenomenon, we simulate a third-order tensor ($N=3$) with a fixed mode size of $n_k=200$ for all $k \in \{1, \dots, 3\}$ and a total of $d=60$ summands. We first generate a true, underlying target tensor $\bTX_{\text{target}}$ with a fixed target rank of $r_{\text{target}} = 5$. To explicitly test the tolerance boundary, the core diagonal of the target tensor is prescribed the values $\{1.0, 0.8, 0.6, 2 \times 10^{-6}, 10^{-6}\}$. This ensures the target possesses an $\mathcal{O}(1)$ norm, while the smallest components intentionally rest exactly at the $\epsilon = 10^{-6}$ relative tolerance boundary. To emulate high-rank, large-amplitude numerical noise that enters and subsequently leaves the domain, we generate $d/2 = 30$ independent noise tensors $\bTX_{\text{noise}}^{(i)}$. Each noise tensor has a fixed rank of $r_{\text{noise}} = 12$, featuring an extreme dynamic range with predefined core diagonal values scaling logarithmically from $10^6$ to $10^{-5}$. The $d$ total summands $\bTX^{(i)}$ are constructed in two symmetric halves to guarantee the mathematical cancellation of these noise components
\begin{align}
    \bTX^{(i)} &= \frac{1}{d} \bTX_{\text{target}} + \bTX_{\text{noise}}^{(i)}, \quad &\text{for } i &= 1, \dots, d/2, \\
    \bTX^{(i+d/2)} &= \frac{1}{d} \bTX_{\text{target}} - \bTX_{\text{noise}}^{(i)}, \quad &\text{for } i &= 1, \dots, d/2.
\end{align}
This construction ensures that the exact mathematical sum $\sum_{i=1}^d \bTX^{(i)}$ perfectly recovers $\bTX_{\text{target}}$ at rank $5$. However, at the exact midpoint of the sequential summation ($i = d/2$), the accumulated, non-cancelled noise forces the nominal rank of the intermediate tensor to swell. Crucially, because the intermediate norm temporarily reaches $\mathcal{O}(10^6)$ due to the noise, applying a relative truncation tolerance of $\epsilon = 10^{-6}$ dynamically forces the absolute truncation threshold up to $\mathcal{O}(1)$. Consequently, the eager rounding  procedure used with \texttt{TuckerAxby} inadvertently truncates and permanently destroys all underlying physical target signals resting below $1$ (namely, the $0.8$, $0.6$, $2 \times 10^{-6}$, and $10^{-6}$ components) long before the canceling noise can conceptually arrive in the second half of the summation. 

We compare four distinct approaches: Eager Rounding (truncating via \texttt{TuckerRounding} immediately after every individual \texttt{TuckerAxby}), Lazy Rounding (forming the naive sum of all $d$ components prior to a single, final truncation), \texttt{KronSum-Tucker}, and \texttt{KRPSum-Tucker} with $p=2$. Relative errors are evaluated explicitly by comparing the full, dense reconstructed tensors. 

\begin{figure}[t]
    \centering
    \begin{subfigure}{0.49\textwidth}
        \centering
  \tikzexternalenable%
  \tikzsetnextfilename{synthetic_cancellation_ranks}%
  \begin{tikzpicture}[font = \plotfontsize\normalfont]
  \pgfplotstableread{graphics/data/Cancellation_Ranks.dat}\tableCancellationRanks

  \begin{axis}[
    width              = \linewidth,
    height             = 5cm,
    enlarge x limits   = 0.02,
    xmin               = 0,
    xmax               = 65,
    ymin               = 0,
    ymax               = 200,
    xtick              = {0, 10, 20, 30, 40, 50, 60},
    xticklabels        = {0, 10, 20, 30, 40, 50, 60},
    ytick              = {0, 50, 100, 150, 200, 250, 300, 350, 400},
    xlabel             = {Summation Step},
    ylabel             = {Max Tucker Rank},
    ylabel style       = {yshift = -.3em},
    ymajorgrids        = true,
    scaled x ticks     = false,
  ]
    \addplot[Naive_Style_line,mark=triangle, mark options={black}] 
    table[
        x=Step, 
        y=RankEager
    ] {\tableCancellationRanks};
  \end{axis}
\end{tikzpicture}%
  \tikzexternaldisable%

        \caption{Eager Rounding Rank Swelling}
        \label{fig:Eager_Rounding_Rank_Swelling}
    \end{subfigure}
    \hfill
        \begin{subfigure}{0.49\textwidth}
        \centering       %
  \tikzexternalenable%
  \tikzsetnextfilename{synthetic_cancellation_times}%
  \begin{tikzpicture}[font = \plotfontsize\normalfont]

\pgfplotstableread{graphics/data/Cancellation_Times.dat}\tableCancellationTimes

\begin{axis}[
    ybar,
    width                 =\linewidth,
    height                = 5cm,
    enlarge x limits      = 0.15,      
    bar width             = 25pt,             
    xtick                 = {0, 1, 2, 3},         
    xticklabels           = {Eager, Lazy, Kron, KRP}, 
    ylabel                = {Total Runtime (s)},
    ymode                 = log,    
    ytick                 = {0.1, 1, 10, 100},
    log origin            = infty,
    ymin                  = 0.1,                    
    ymax                  = 100,                  
    ymajorgrids           = true,
    grid style            = {dotted, gray!60}
]
    \addplot[Naive_Style_Bar, bar shift=0pt, mark=triangle, mark options={black, xshift=0pt, yshift=-8pt}] 
        table[x expr=\coordindex, y=Time, restrict expr to domain={\coordindex}{0:0}] {\tableCancellationTimes};
        
    \addplot[Naive_Style_Bar, bar shift=0pt,mark=triangle, mark options={black, xshift=0pt, yshift=-8pt}] 
        table[x expr=\coordindex, y=Time, restrict expr to domain={\coordindex}{1:1}] {\tableCancellationTimes};
        
    \addplot[Kron_Style_Bar, bar shift=0pt,mark=square*, mark options={black, xshift=0pt, yshift=-8pt}] 
        table[x expr=\coordindex, y=Time, restrict expr to domain={\coordindex}{2:2}] {\tableCancellationTimes};
        
    \addplot[KRP_Style_Bar, bar shift=0pt,mark=*, mark options={black, xshift=0pt, yshift=-8pt}] 
        table[x expr=\coordindex, y=Time, restrict expr to domain={\coordindex}{3:3}] {\tableCancellationTimes};

\end{axis}
\end{tikzpicture}%
  \tikzexternaldisable%

        \caption{Runtimes}        \label{fig:Synthetic_cancellation_times}
    \end{subfigure}   
    \centering   
    \begin{subfigure}{0.49\textwidth}
        \centering 
  \tikzexternalenable%
  \tikzsetnextfilename{synthetic_cancellation_errors}%
  \begin{tikzpicture}[font = \plotfontsize\normalfont]
        \node[inner sep=0pt] {\pgfplotstabletypeset[
            col sep=space,
            every row/.style={after row=\hline},
            every last row/.style={after row=\hline},
            every head row/.style={
                before row={
                    \hline
                },
                after row=\hline
            },
            columns={Method, Error},
            columns/Method/.style={
                column name={Method}, 
                string type,
                column type=|c| 
            },            
            columns/Error/.style={
                column name={Relative Error}, 
                sci, sci zerofill, precision=2,
                column type=c|
            }
        ]{graphics/data/Cancellation_Errors.dat}};
\end{tikzpicture}%
  \tikzexternaldisable%

        \caption{Relative Approximation Error}        \label{fig:Synthetic_cancellation_errors}
    \end{subfigure}
      \vspace{.5\baselineskip}
  \tikzexternalenable%
  \tikzsetnextfilename{Legend}%
  \begin{tikzpicture}[font = \plotfontsize\normalfont]
  \begin{axis}[%
    hide axis, 
    scale only axis,
    width          = .8\linewidth,
    height         = .2\linewidth,
    xmin           = 0,
    xmax           = 1,
    ymin           = 0,
    ymax           = 1,
    legend columns = 3,    
    legend style   = {
      at     = {(0.5,0.5)}, 
      anchor = center,
      draw   = none,       
      /tikz/every even column/.append style = {column sep = 0.5cm}
    }
  ]
    
    
    \addlegendimage{Kron_Style_line, mark=square*, mark options={black}}
    \addlegendentry{\texttt{KronSum-Tucker}}

    \addlegendimage{KRP_Style_line,mark=*, mark options={black}}
    \addlegendentry{\texttt{KRPSum-Tucker}}

    \addlegendimage{Naive_Style_line,mark=triangle, mark options={black}}
    \addlegendentry{\texttt{TuckerAxby}+\texttt{TuckerRounding}}

  \end{axis}
\end{tikzpicture}%
  \tikzexternaldisable%

      \caption{Performance comparison of sequential truncation versus holistic summation: 
Eager deterministic rounding suffers from severe intermediate rank swelling and accuracy loss (relative error of $6.72 \cdot 10^{0}$) due to delayed mathematical cancellation. 
By contrast, the randomized sketching approaches inherently bypass this intermediate growth with a low error relative to the true exact tensor ($\sim 2 \cdot 10^{-11}$) while running an order of magnitude faster than the eager baseline and maintaining a distinct computational advantage over the lazy baseline.}
\label{fig:Synthetic_cancellation}
  
\end{figure}
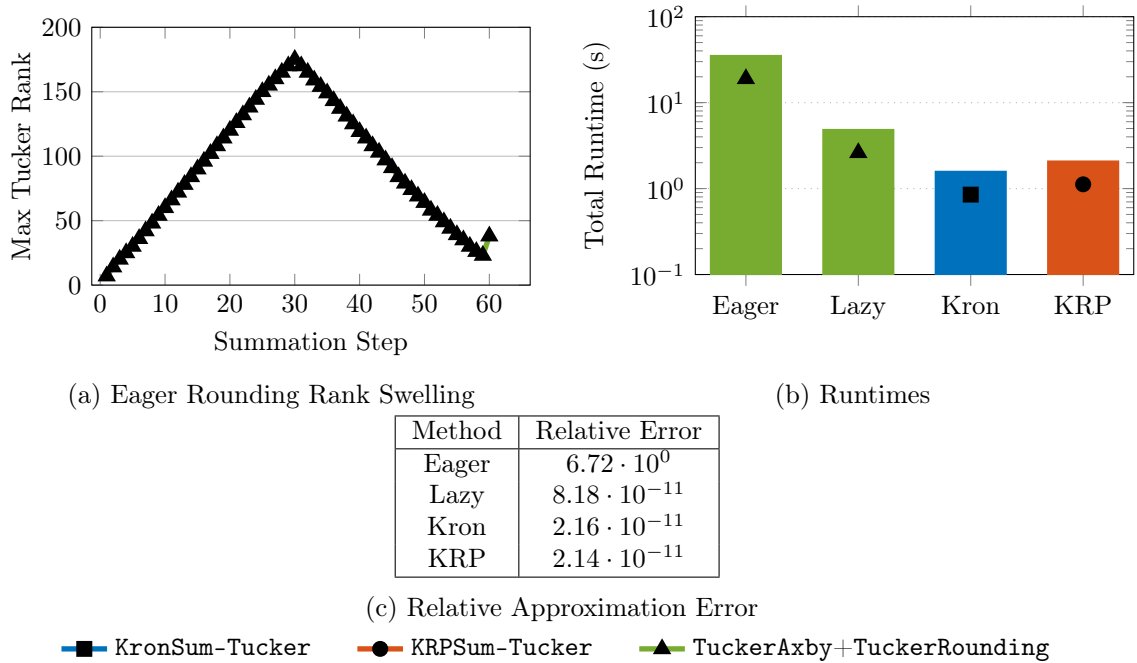

As \Cref{fig:Synthetic_cancellation} demonstrates, deterministic eager rounding fails catastrophically. Intermediate norm inflation artificially raises the absolute truncation threshold, permanently destroying delicate $\mathcal{O}(1)$ target signals before the adversarial noise cancels, while simultaneously imposing severe computational penalties from rank swelling. While lazy rounding evaluates the sum holistically to preserve accuracy, it still requires forming a massively swollen intermediate core. The proposed randomized Kronecker and KRP sketching strategies prove superior by jointly solving both vulnerabilities: they perfectly preserve the underlying signal and entirely bypass the computational bottlenecks of intermediate rank swelling. As a result we will set lazy rounding as our baseline deterministic method to avoid rank swelling and potential accuracy issues.

\subsection{Cookie Problem}
\label{subsec:cookie_problem}

For our first large-scale experiment, we consider the \emph{cookie problem}~\cite{KreT11, Tob12}, which has been used as a benchmark for randomized rounding algorithms in the tensor-train format~\cite{AldBCetal23}. We adapt their implementation directly to the Tucker format in order to compare our randomized approach with standard Tucker rounding.
We consider the stationary heat equation posed on the unit square domain $\Omega = [0,L]^2$, with diffusion coefficient $\sigma(x,\bxi)$ depending on a set of independent parameters $\bxi = (\xi_1,\dots,\xi_P) \in \Gamma \subset \R^P$. The goal is to determine the solution $u(x,\bxi)$ satisfying
\begin{align}
    -\nabla \cdot (\sigma(x, \bxi) \nabla u(x, \bxi)) &= f(x) \quad \text{in } \Omega, \\
    u(x, \bxi) &= 0 \quad \text{on } \partial\Omega,
\end{align}
where $f(x)$ represents a constant source term. The domain $\Omega$ contains $P$ mutually disjoint circular subdomains $\mathcal{D}_1, \dots, \mathcal{D}_P$, which correspond to the geometric inclusions shown in the experiment setup. The diffusion coefficient is modeled as piecewise constant, taking a base value of $1$ outside the inclusions and a parameter-dependent value $1 + \xi_\mu$ inside the $\mu$-th inclusion $\mathcal{D}_\mu$
$$
\sigma(x, \bxi) = 1 + \sum_{\mu=1}^P \xi_\mu \chi_{\mathcal{D}_\mu}(x),
$$
where $\chi_{\mathcal{D}_\mu}$ denotes the characteristic function of the subdomain $\mathcal{D}_\mu$.

\begin{figure}[t]
    \centering
  \tikzexternalenable%
  \tikzsetnextfilename{Cookie_Domain}%
  \begin{tikzpicture}[scale=2.2]
    \definecolor{slate}{RGB}{38, 50, 56}      %
    \definecolor{amber}{RGB}{255, 183, 77}
    \definecolor{highlight}{RGB}{255, 255, 255}

    \filldraw[fill=slate, draw=black, thick, rounded corners=2pt] (0,0) rectangle (3,3);

    \foreach \x/\y/\mu in {0.75/0.75/1, 2.25/0.75/2, 0.75/2.25/3, 2.25/2.25/4} {
        \filldraw[fill=amber, draw=white, line width=0.8pt] (\x,\y) circle (0.42);
        
        \node[text=slate, font=\boldmath\small] at (\x,\y) {$D_{\mu}$};
        \node[text=slate, font=\scriptsize] at (\x, \y - 0.22) {$1 + \xi_{\mu}$};
    }

    \node[draw=white, fill=white, fill opacity=0.95, 
          text=slate, font=\bfseries\small, 
          rounded corners=3pt, inner sep=8pt, line width=1.5pt] 
        at (1.5, 1.5) {$\Omega: \sigma(\mathbf{x}, \boldsymbol{\xi}) = 1$};

    \node[below, text=slate, font=\footnotesize] at (0, -0.05) {$0$};
    \node[below, text=slate, font=\footnotesize] at (3, -0.05) {$L$};
    \node[left,  text=slate, font=\footnotesize] at (-0.05, 0) {$0$};
    \node[left,  text=slate, font=\footnotesize] at (-0.05, 3) {$L$};

\end{tikzpicture}%
  \tikzexternaldisable%

    \caption{Schematic of the spatial domain $\Omega = [0,L]^2$ for the 4-parameter \emph{cookie problem}. The domain features a background region where the diffusion coefficient is constant ($\sigma = 1$), and $P=4$ mutually disjoint circular subdomains arranged in a grid. Inside each subdomain $\mathcal{D}_\mu$, the diffusion coefficient is modeled as a parameter-dependent constant $\sigma(x, \boldsymbol{\xi}) = 1 + \xi_\mu$. This geometry gives rise to the spatial stiffness matrices $\mathbf{A}_0, \dots, \mathbf{A}_3$ utilized in the tensor-structured linear system.}
    \label{fig:cookie_domain}
\end{figure}
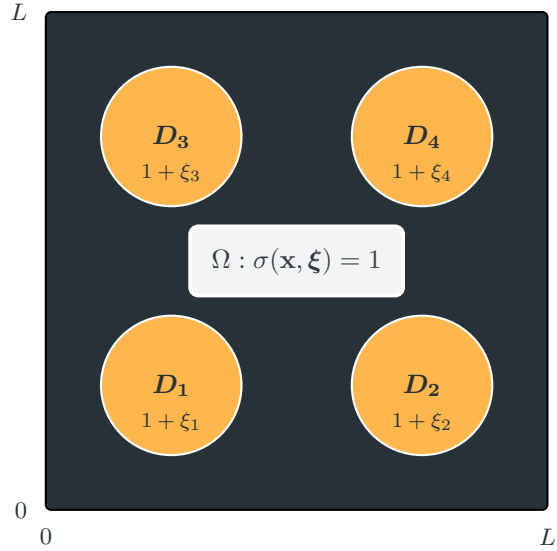

Discretizing the spatial domain using finite elements with $N_x$ basis functions $\{\phi_i\}_{i=1}^{N_x}$, we approximate the solution as $u_h(x, \bxi) = \sum_{i=1}^{N_x} u_i(\bxi) \phi_i(x)$. The weak formulation leads to a parameter-dependent linear system of the form
\begin{equation}
    \left( \bA_0 + \sum_{\mu=1}^P \xi_\mu \bA_\mu \right) \bu(\bxi) = \bb,
    \label{eq:param_system}
\end{equation}
where $\bu(\bxi) \in \R^{N_x}$ is the vector of coefficients. The matrices $\bA_0, \bA_\mu \in \R^{N_x \times N_x}$ and the vector $\bb \in \R^{N_x}$ are defined by the integrals
$$
(\bA_0)_{ij} = \int_\Omega \nabla \phi_i \cdot \nabla \phi_j \, dx, \quad (\bA_\mu)_{ij} = \int_{\mathcal{D}_\mu} \nabla \phi_i \cdot \nabla \phi_j \, dx, \quad (\bb)_i = \int_\Omega f \phi_i \, dx.
$$
In our experiments, we utilize pre-computed stiffness matrices $\bA_0, \dots, \bA_P$ and load vector $\bb$. To solve this system for a range of parameter values simultaneously, we employ a tensor product approach. We discretize the parameter space by sampling each parameter $\xi_\mu$ at $N$ points, denoted by the set $\{\xi_{\mu}^{(k)}\}_{k=1}^N$. This sampling forms a grid of $N^P$ total configurations. The discrete unknowns are collected into a $(P+1)$-dimensional tensor $\bTX \in \R^{N_x \times N \times \dots \times N}$. The parameter-dependent linear system from~\cref{eq:param_system} can be formulated globally as a linear operator $\mathscr{L}$ acting on this solution tensor. By defining $\bD = \text{diag}(\xi^{(1)}, \dots, \xi^{(N)}) \in \R^{N \times N}$ as the diagonal matrix containing the parameter samples and $\bI \in \R^{N \times N}$ as the identity matrix, the operator $\mathscr{L}$ can be expressed algebraically as a sum of Kronecker products. For example, for the specific case of $P=2$ parameters, the global operator takes the form
\begin{equation}\mathscr{L} = \bA_0 \otimes \bI \otimes \bI + \bA_1 \otimes \bD \otimes \bI + \bA_2 \otimes \bI \otimes \bD.\label{eq:tensor_operator}
\end{equation}
When applied to $\bTX$, the individual matrices within each Kronecker term act along their corresponding tensor modes. The right-hand side of the system is the $(P+1)$-dimensional tensor $\bTB$. It is constructed via the outer product (denoted by $\circ$) of the spatial load vector $\bb \in \R^{N_x}$ and the parameter domain vectors. This yields a rank-1 tensor representing the constant source term distributed across the entire parameter grid
$$\bTB = \bb \circ \mathbf{1} \circ \mathbf{1},$$
where $\mathbf{1} \in \R^N$ is a vector of ones. We solve the resulting tensor-structured linear system, formally evaluated as $\mathscr{L}(\bTX) = \bTB$, using the \texttt{Tucker-GMRES} solver described in \Cref{alg:tucker_gmres}.

\begin{algorithm}[t]
    \caption{Right-Preconditioned \texttt{Tucker-GMRES} with \texttt{Tucker-Rounding}}
    \label{alg:tucker_gmres}
    \SetKwInOut{Input}{Input}
    \SetKwInOut{Output}{Output}
    \DontPrintSemicolon

    \Input{Linear operator $\bTA$, RHS $\bTB$, Preconditioner $\bTM^{-1}$, Tolerance $\tau$, Max Iterations $K$, Truncation Operator $\texttt{RoundSum}(\cdot)$}
    \Output{Approximate solution $\bTX$}

    $\bTR_0 \gets \bTB$\;
    $\beta \gets \|\bTR_0\|_F$\;
    $\bTV_1 \gets \bTR_0 / \beta$\;
    Initialize Hessenberg matrix $\bH \in \R^{(K+1) \times K}$\;

    \For{$k \gets 1$ \KwTo $K$}{
        $\bTW \gets \bTM^{-1}(\bTV_k)$ \tcp*{Apply preconditioner}
        $\bTW \gets \bTA(\bTW)$ \tcp*{Apply operator}
        
        $\bTW \gets \texttt{RoundSum}(\{\bTW\}, \{1\})$ \tcp*{Optional compression}

        \For{$j \gets 1$ \KwTo $k$}{
            $h_{j,k} \gets \langle \bTV_j, \bTW \rangle$ \tcp*{Tensor inner product}
        }
        
        $\bTW \gets \texttt{RoundSum}(\{\bTV_1, \dots, \bTV_k, \bTW\}, \{-h_{1,k}, \dots, -h_{k,k}, 1\})$\;
        
        $h_{k+1, k} \gets \|\bTW\|_F$\;
        
        \If{$h_{k+1, k} \approx 0$}{
             Break\;
        }
        
        $\bTV_{k+1} \gets \bTW / h_{k+1, k}$\;
        
        $\by_k \gets \arg\min_{\by} \| \beta \be_1 - \bH_{1:k+1, 1:k} \by \|_2$\;
        
        Compute residual error $\varepsilon = \| \beta \be_1 - \bH_{1:k+1, 1:k} \by_k \|_2$\;
        \If{$\varepsilon < \tau \|\bTB\|_F$}{
            Break\;
        }
    }
    
    $\bTY \gets \texttt{RoundSum}(\{\bTV_1, \dots, \bTV_k\}, \by_k)$\;
    
    $\bTX \gets \bTM^{-1}(\bTY)$\;
    
    $\bTX \gets \text{\texttt{TuckerRounding}}(\bTX, \tau)$\;

    \Return $\bTX$\;
\end{algorithm}

In our numerical experiments, the parameters $\xi_\mu$ are sampled linearly within the interval $[1, 10]$, and the number of parameter samples is systematically varied across $N \in \{8, 16, 32, 64, 128\}$. For preconditioning, we employ a reference-based global operator $\mathscr{M}^{-1}$ acting on the spatial mode, which approximates the system at the fixed configuration $\xi_\mu = 1$. This is realized by first forming the reference spatial matrix $\bM = \bA_{0} + \sum_{\mu=1}^{P} \bA_{\mu}$ and computing its LU factorization. The global preconditioner is then applied to a given tensor $\bTY$ via the mode-1 tensor-matrix product $\mathscr{M}^{-1}(\bTY) = \bTY \times_1 \bM^{-1}$, utilizing the precomputed LU factors to avoid explicitly forming the inverse. For the iterative solver, the \texttt{Tucker-GMRES} convergence tolerance is set to $10^{-5}$, with a maximum allowance of 20 iterations. Crucially, to prevent accuracy loss during the intermediate stages, the internal tensor truncation tolerance applied during the Arnoldi summation steps is set tighter than the convergence criteria, specifically to $10^{-7}$. For the randomized algorithms, an oversampling parameter of $p = 5$ is used. To ensure reliable runtime comparisons and mitigate stochastic variance, all numerical experiments are averaged over 3 independent trials.

As illustrated in \Cref{fig:GMRES_example}, while \texttt{KRPSum-Tucker} demonstrates robust and predictable performance, \texttt{KronSum-Tucker} exhibits degraded efficiency. This performance discrepancy is not an inherent flaw of Kronecker sketching itself, but rather a direct consequence of the highly skewed dimensions characteristic of this specific PDE application, where the spatial dimension $N_x$ is vastly larger than the parameter sample dimensions $N$. In such heavily skewed regimes, the standard \textit{factor reuse strategy} becomes algorithmically susceptible and loses its computational advantages. We emphasize that these results should not be interpreted as a general discouragement of Kronecker sketching; instead, they highlight a critical structural caveat: factor reuse is inappropriate for tensors with severely skewed mode dimensions. To properly leverage Kronecker sketching in this parametric PDE setting, one must forgo factor reuse and instead adopt different sketching heuristics~\cite{MinLB24}. Consequently, for the highly skewed geometries evaluated in this section, the KRP formulation provides a much more natural and robust randomized alternative.

\begin{figure}[p]
    \centering
    \begin{subfigure}{0.99\textwidth}
        \centering
  \tikzexternalenable%
  \tikzsetnextfilename{GMRES_Runtime}%
  \begin{tikzpicture}[font = \plotfontsize\normalfont]

\pgfplotstableread{graphics/data/GMRES_Runtime.dat}\tableRuntime

\begin{axis}[
    ybar,
    width                   = \linewidth,
    height                  = 5cm,
    enlarge x limits        = 0.15,
    bar width               = 20pt,
    xtick                   = data,
    xticklabels from table  = {\tableRuntime}{N},
    xlabel                  = {Number of samples $N$},
    x tick label style      = {/pgf/number format/1000 sep=},
    ylabel                  = {Total Runtime (s)},
    ytick                   = {0.1, 1, 10, 100, 1000, 10000},
    ymode                   = log,      
    log origin              = infty,
    ymin                    = 0.1,      
    ymajorgrids             = true,
]
    \addplot[Kron_Style_Bar, mark=square*, mark options={black,xshift=-22pt, yshift=-8pt}] table[x expr=\coordindex, y=KronSum] {\tableRuntime};

    \addplot[KRP_Style_Bar,mark=*, mark options={black,xshift=0pt, yshift=-8pt}] table[x expr=\coordindex, y=KRPSum] {\tableRuntime};

    \addplot[Naive_Style_Bar,mark=triangle, mark options={black,xshift=22pt, yshift=-8pt}] table[x expr=\coordindex, y=NaiveSum] {\tableRuntime};

\end{axis}
\end{tikzpicture}%
  \tikzexternaldisable%

        \caption{Runtime}
        \label{fig:GMRES_Runtime}
    \end{subfigure}
    \vspace{0.5cm}

        \begin{subfigure}{0.99\textwidth}
        \centering
  \tikzexternalenable%
  \tikzsetnextfilename{GMRES_Ranks}%
  \definecolor{matDet}{HTML}{77AC30}
\definecolor{matKron}{HTML}{0072BD}
\definecolor{matKRP}{HTML}{D95319}

\begin{tikzpicture}[font = \plotfontsize\normalfont]
  \pgfplotstableread{graphics/data/GMRES_Ranks.dat}\tableRanks

  \begin{axis}[
    width              = 0.99\linewidth,
    height             = 5cm,
    enlarge x limits   = 0.025,
    xmin               = 0,
    xmax               = 20,
    ymin               = 0,
    ymax               = 140,
    xtick              = {0, 5, 10, 15, 20},
    ytick              = {0, 20, 40, 60, 80, 100,120,140},
    xlabel             = {Iteration},
    ylabel             = {Tucker Rank},
    ylabel style       = {yshift = -.3em},
    ymajorgrids        = true,
    legend pos         = north west,
    legend style       = {font=\tiny, row sep=0.5em},
    scaled x ticks     = false,
  ]
    
    \addplot[matDet, only marks, mark=diamond,  mark size=5.0pt, mark options={line width=0.8pt}, forget plot] table[x=Iteration, y=Naive_D1] {\tableRanks};
    \addplot[matDet, only marks, mark=square,   mark size=5.0pt, mark options={line width=0.8pt}, forget plot] table[x=Iteration, y=Naive_D2] {\tableRanks};
    \addplot[matDet, only marks, mark=triangle, mark size=5.0pt, mark options={line width=0.8pt}, forget plot] table[x=Iteration, y=Naive_D3] {\tableRanks};
    \addplot[matDet, only marks, mark=o,        mark size=5.0pt, mark options={line width=0.8pt}, forget plot] table[x=Iteration, y=Naive_D4] {\tableRanks};

    \addplot[matKron, only marks, mark=diamond,  mark size=4.0pt, mark options={line width=0.8pt}, forget plot] table[x=Iteration, y=Kron_D1] {\tableRanks};
    \addplot[matKron, only marks, mark=square,   mark size=4.0pt, mark options={line width=0.8pt}, forget plot] table[x=Iteration, y=Kron_D2] {\tableRanks};
    \addplot[matKron, only marks, mark=triangle, mark size=4.0pt, mark options={line width=0.8pt}, forget plot] table[x=Iteration, y=Kron_D3] {\tableRanks};
    \addplot[matKron, only marks, mark=o,        mark size=4.0pt, mark options={line width=0.8pt}, forget plot] table[x=Iteration, y=Kron_D4] {\tableRanks};

    \addplot[matKRP, only marks, mark=diamond,  mark size=4.0pt, mark options={line width=0.8pt}, forget plot] table[x=Iteration, y=KRP_D1] {\tableRanks};
    \addplot[matKRP, only marks, mark=square,   mark size=4.0pt, mark options={line width=0.8pt}, forget plot] table[x=Iteration, y=KRP_D2] {\tableRanks};
    \addplot[matKRP, only marks, mark=triangle, mark size=4.0pt, mark options={line width=0.8pt}, forget plot] table[x=Iteration, y=KRP_D3] {\tableRanks};
    \addplot[matKRP, only marks, mark=o,        mark size=4.0pt, mark options={line width=0.8pt}, forget plot] table[x=Iteration, y=KRP_D4] {\tableRanks};

    \addlegendimage{only marks, mark=diamond, mark options={color=black, fill=black}, mark size=4pt}
    \addlegendentry{Mode 1}

    \addlegendimage{only marks, mark=square, mark options={color=black}, mark size=4pt}
    \addlegendentry{Mode 2}

    \addlegendimage{only marks, mark=triangle, mark options={color=black}, mark size=4pt}
    \addlegendentry{Mode 3}

    \addlegendimage{only marks, mark=o, mark options={color=black}, mark size=4pt} 
    \addlegendentry{Mode 4}

  \end{axis}
\end{tikzpicture}%
  \tikzexternaldisable%

        \caption{Max Tucker Rank}
        \label{fig:GMRES_Ranks}
    \end{subfigure}

\vspace{0.5cm}
    
    \begin{subfigure}{0.49\textwidth}
        \centering
  \tikzexternalenable%
  \tikzsetnextfilename{GMRES_Speedup}%
  \begin{tikzpicture}[font = \plotfontsize\normalfont]
  \pgfplotstableread{graphics/data/GMRES_Speedup.dat}\tableSpeedup

  \begin{axis}[
    width              = \linewidth,
    height             = 4cm,
    enlarge x limits   = 0.15,
    xmin               = 8,
    xmax               = 128,
    ymin               = 0,
    ymax               = 12,
    xmode              = log,
    xtick              = {8, 16, 32, 64, 128},
    xticklabels        = {8, 16, 32, 64, 128},
    ytick              = {0, 2, 4, 6, 8, 10, 12},
    xlabel             = {$N$},
    ylabel             = {Speedup},
    ylabel style       = {yshift = -.3em},
    ymajorgrids        = true,
    scaled x ticks     = false,
  ]
    
    \addplot[Kron_Style_line,mark=square*, mark options={black}] 
    table[
        x=N, 
        y=KronSpeedup, 
        y error minus=KronErrMin, 
        y error plus=KronErrMax
    ] {\tableSpeedup};

    \addplot[KRP_Style_line,mark=*, mark options={black}] 
    table[
        x=N, 
        y=KRPSpeedup, 
        y error minus=KRPErrMin, 
        y error plus=KRPErrMax
    ] {\tableSpeedup};

    \addplot[Naive_Style_line, mark=triangle*, mark options={black}] coordinates {(8,1) (16,1) (32,1) (64,1) (128,1)};

  \end{axis}
\end{tikzpicture}%
  \tikzexternaldisable%

        \caption{Speed up}
        \label{fig:GMRES_Speedup}
    \end{subfigure}
    \hfill
        \begin{subfigure}{0.49\textwidth}
        \centering
  \tikzexternalenable%
  \tikzsetnextfilename{GMRES_Accuracy}%
  \begin{tikzpicture}[font = \plotfontsize\normalfont]
  \pgfplotstableread{graphics/data/GMRES_Accuracy.dat}\tableAccuracy

  \begin{axis}[
    width              = \linewidth,
    height             = 4cm,
    enlarge x limits   = 0.15,
    xmin               = 0,
    xmax               = 140,
    ymode              = log,
    ymin               = 1e-6, 
    ymax               = 1e-1,
    xmode              = log,
    xtick              = {8, 16, 32, 64, 128},
    xticklabels        = {8, 16, 32, 64, 128},
    xlabel             = {$N$},
    ylabel             = {Relative Error}, 
    ylabel style       = {yshift = -.3em},
    ymajorgrids        = true,
    scaled x ticks     = false,
  ]

    \addplot[KRP_Style_line, mark=*, mark options={black}] 
    table[
        x=N, 
        y=KRPMean, 
        y error minus=KRPMin, 
        y error plus=KRPMax
    ] {\tableAccuracy};

    \addplot[Kron_Style_line,mark=square*, mark options={black}] 
    table[
        x=N, 
        y=KronMean, 
        y error minus=KronMin, 
        y error plus=KronMax
    ] {\tableAccuracy};

    \addplot[Naive_Style_line, dashed, mark=triangle*,mark size=2.5pt, mark options={solid,black}] 
    table[
        x=N, 
        y=NaiveMean, 
        y error minus=NaiveMin, 
        y error plus=NaiveMax
    ] {\tableAccuracy};

  \end{axis}
\end{tikzpicture}%
  \tikzexternaldisable%

        \caption{Accuracy}
        \label{fig:GMRES_Accuracy}
    \end{subfigure}
  \tikzexternalenable%
  \tikzsetnextfilename{Legend}%
  \begin{tikzpicture}[font = \plotfontsize\normalfont]
  \begin{axis}[%
    hide axis, 
    scale only axis,
    width          = .8\linewidth,
    height         = .2\linewidth,
    xmin           = 0,
    xmax           = 1,
    ymin           = 0,
    ymax           = 1,
    legend columns = 3,    
    legend style   = {
      at     = {(0.5,0.5)}, 
      anchor = center,
      draw   = none,       
      /tikz/every even column/.append style = {column sep = 0.5cm}
    }
  ]
    
    
    \addlegendimage{Kron_Style_line, mark=square*, mark options={black}}
    \addlegendentry{\texttt{KronSum-Tucker}}

    \addlegendimage{KRP_Style_line,mark=*, mark options={black}}
    \addlegendentry{\texttt{KRPSum-Tucker}}

    \addlegendimage{Naive_Style_line,mark=triangle, mark options={black}}
    \addlegendentry{\texttt{TuckerAxby}+\texttt{TuckerRounding}}

  \end{axis}
\end{tikzpicture}%
  \tikzexternaldisable%

        \caption{Performance evaluation of the randomized \texttt{Tucker-GMRES} solver for the parameter-dependent \emph{cookie problem}: 
    \texttt{KRPSum-Tucker} demonstrates superior efficiency and robustness compared to both \texttt{KronSum-Tucker} and \texttt{TuckerAxby}+\texttt{TuckerRounding}. 
    Across varying parameter sample sizes ($N$), the \texttt{KRPSum-Tucker} achieves up to a $11\times$ computational speedup while preserving high accuracy with relative errors resting near $10^{-5}$. 
    In contrast, \texttt{KronSum-Tucker} struggles in this specific setting, yielding higher relative errors ($\sim 10^{-2}$) and minimal runtime improvements over the standard rounding baseline due to a maximum 30 iteration cap. \Cref{fig:GMRES_Ranks} demonstrates that the highly skewed mode in the Kronecker case has inflated rank compared to \texttt{KRPSum-Tucker} and \texttt{TuckerAxby}+\texttt{TuckerRounding}.}

      \vspace{.5\baselineskip}
      \label{fig:GMRES_example}
\end{figure}

\subsection{Linear conservation law}
\label{subsec:Linear_conservation_law}
We consider a parametric 1D linear conservation law
\begin{equation}
    \frac{\partial f}{\partial t}+\xi_x\frac{\partial f}{\partial x}=0,\qquad x\in[0,2\pi],\ \bxi\in[-\Xi,\Xi]^3,
    \label{eq:transport_equation}
\end{equation}
where $f:=f(x,t;\bxi)$ and $\bxi=(\xi_x,\xi_y,\xi_z)\in\R^3$, we also assume that $f$ decays on the parameter $\bxi$ as $|\bxi|\rightarrow\infty$. We enforce periodic boundary conditions in physical space $x$ and assume a large enough numerical domain ($\Xi$) for parameter $\bxi$ to allow for sufficient decay of the solution. High-dimensional transport equations similar to \Cref{eq:transport_equation} appear in many applications, such as rarefied gas dynamics, neutron transport in high energy density systems, and plasma modeling. In many rank-adaptive methods~\cite{EinKKetal25}, the distribution function is stored using low-rank tensor decompositions to exploit low-rank structures of the solution in $\bxi$. In particular, the Tucker decomposition has been recently used to solve 3D conservation laws~\cite{AppC25, ElKQCetal25, NakCE25}. Assuming an initial condition $f_0:=f_0(x,\bxi)$, the method of characteristics gives the exact solution to \Cref{eq:transport_equation},
\begin{equation}
    f(x,t;\bxi)=f_0(x-t\xi_x).
    \label{eq:exact_solution_transport}
\end{equation}

Following the work presented in~\cite{GalNPetal25} for the 2D parameter space case, we use a Nodal Discontinuous Galerkin (NDG) method to discretize in $x$~\cite{HesW08}. At each spatial node we store a Tucker decomposition of the solution in $\bxi$ that evolves in time. Classical time-stepping methods (e.g., Runge-Kutta, linear multistep) used to solve \Cref{eq:transport_equation} require the summation of the transport term and/or solution evaluated at many times/stages. Higher-order time integrators require a greater number of tensor summations, posing a significant computational bottleneck. The efficient sketching-based method presented in this paper offers a way to mitigate this bottleneck, allowing one to preserve the efficiency that comes from the NDG and low-rank Tucker frameworks.

For brevity, we only state the key equations and modifications made to the NDG method in \cite{GalNPetal25} and refer the reader to the original paper for more details. We then show the numerical performance resulting from using our proposed sketching-based method for the tensor summations. Discretizing the spatial domain into $N_x$ intervals,
\[0=x_{\frac{1}{2}}<x_{\frac{3}{2}}<...<x_{N_x-\frac{1}{2}}<x_{N_x+\frac{1}{2}}=2\pi,\]
generates a mesh consisting of elements $I_i\coloneqq[x_{i-\frac{1}{2}},x_{i+\frac{1}{2}}]$, $i=1,2,...,N_x$. For simplicity, we assume a uniform mesh with $h_x=x_{i+\frac{1}{2}}-x_{i-\frac{1}{2}}$ for all $i$. Define the finite dimensional discrete space
\begin{equation}
    V_h^k\coloneqq \{g\in L^2(\Omega_x)\ :\  g|_{I_i}\in \mathbb{P}^k(I_i), 1\leq i\leq N_x\},
\end{equation}
where $k\geq 0$ is a non-negative integer, and $\mathbb{P}^k(I_i)$ denotes the space of polynomials of degree at most $k$ on the element $I_i$. For each element $I_i$ ($1\leq i\leq N_x$), we consider the associated polynomial basis consisting of Lagrange polynomials $L^i_q$ locally defined over $k+1$ Gauss-Legendre quadrature nodes $\{x^i_s\in I_i\ :\ s=1,...,k+1\}$. We desire the weak solution $f_h(\cdot,t;\bxi)\in V_h^k$ such that for all $i=1,2,...,N_x$ and for all $\varphi\in V_h^k$,
\begin{equation}\label{eq:weakform}
    \int_{I_i}{\frac{\partial f_h}{\partial t}\varphi dx} + \widehat{(\xi_{x}f_h)}_{i+\frac{1}{2}}\varphi^-_{i+\frac{1}{2}} - \widehat{(\xi_{x}f_h)}_{i-\frac{1}{2}}\varphi^+_{i-\frac{1}{2}} - \int_{I_i}{\xi_{x}f_h\frac{d\varphi}{dx}dx}= 0,
\end{equation}
where subscript $i\pm\frac{1}{2}$ indicates an approximation as $x\rightarrow x_{i\pm\frac{1}{2}}$, superscript $\pm$ indicates the right or left limit, and $\widehat{\xi_{x}f_h}$ are the upwind numerical fluxes. Decomposing the weak solution (restricted to an element $I_i$) using the associated polynomial basis,
\begin{equation}\label{eq:fimn}
    f_h(x\in I_i,t;\bxi) = \sum\limits_{q=1}^{k+1}{C^{i}_q(\bxi,t)L_q^i(x)}.
\end{equation}

Substituting \Cref{eq:fimn} into \Cref{eq:weakform}, integrating against each basis function $L^i_p\in V_h^k$, and approximating all integrals using a Gaussian quadrature yields a differential equation for the coefficients of the form
\begin{equation}\label{eq:Cipmn}
\frac{dC^i_p}{dt} = \sum\limits_{\alpha=1}^{A(k)}\sum\limits_{\beta\in S_i}{\omega^{\beta}_{\alpha}C^{\beta}_{\alpha}},
\end{equation}
where $A(k)=2(k+1)$ due to both numerical fluxes, $S_i$ is a stencil containing $i$ (e.g., $S_i=\{i-1,i,i+1\}$), and $\omega_{\alpha}^{\beta}$ are some coefficients. When using the Lagrange polynomials as the choice of basis, the coefficients are identical to nodal values of the solution, that is, $f_h(x^i_p\in I_i,\bxi,t)=C^{i}_{p}(\bxi,t)$ for all $i=1,...,N_x$ and all $p=1,...,k+1$, by \Cref{eq:fimn}. That is, computing the coefficients in turn produces the nodal values of the solution.

To discretize in parameter space $\bxi$, we assume $N_{\bxi}$ equally spaced points in $\xi_x,\xi_y,\xi_z$; the width is denoted by $h_v=2\Xi/N_{\bxi}$. For each element $I_i$ ($i=1,...,N_x$) and each node $x_p^i$ ($p=1,...,k+1$), we update a third-order tensor solution $\bC^i_p(t)\in\R^{N_{\bxi}\times N_{\bxi} \times N_{\bxi}}$. We assume each $\bC^i_p(t)$ is decomposed as a third-order Tucker tensor. Integrating \Cref{eq:Cipmn} in time using forward Euler,
\begin{equation}\label{eq:Cipmn_fEuler}
    (\bC^i_p)^{n+1}=(\bC^i_p)^{n}+\Delta t\sum\limits_{\alpha=1}^{A(k)}\sum\limits_{\beta\in S_i}{\omega^{\beta}_{\alpha}(\bC^{\beta}_{\alpha})^{n}}.
\end{equation}

Looking at \Cref{eq:Cipmn_fEuler}, our proposed sketching-based summation technique can be used to significantly improve the computational efficiency of evaluating the right-hand side. One can imagine even greater savings for higher-order time-stepping methods where a greater number of summations is needed.

First, we numerically verify the spatial accuracy of the NDG methods when using (naive) \texttt{TuckerAxby}+\texttt{TuckerRounding}, \texttt{KRPSum-Tucker}, and \texttt{KronSum-Tucker}. We assume a smooth rank-1 initial condition,
\begin{equation}
    f(x,0;\bxi)= \frac{1+0.5\sin{(x)}}{(2\pi)^{3/2}}\exp(-|\bxi|^2/2).
\end{equation}

We set the parameter space mesh $N_{\bxi}=32$ with $\Xi=6$, tolerance $10^{-6}$, final time $T_f=1.0$, and time-step size $\Delta t=\vartheta\Delta x^{k+1}$, where $\vartheta=0.1/((2k+3)\Xi)$ with $k=0,1,2$. The NDG method of (spatial) order 1, 2, 3 are denoted by NDG1, NDG2, NDG3, respectively. As seen in \Cref{fig:NDG_check}, we observe the expected order of convergence (in $L^1$) under spatial mesh refinement. The $L^1$ errors were indistinguishable when using all three techniques: \texttt{KronSum-Tucker}, \texttt{KRPSum-Tucker}, \texttt{TuckerAxby+TuckerRounding}.

\begin{figure}[t!]
    \centering
  \tikzexternalenable%
  \tikzsetnextfilename{NDG_check}%
  \pgfplotstabletypeset[
    col sep=space,
    empty cells with={-},
    every row/.style={after row=\hline},
    every last row/.style={after row=\hline},
    every head row/.style={
        before row={
            \hline
        },
        after row=\hline
    },
    columns={Basis, Nx, ErrDet, OrdDet},
    columns/Basis/.style={
        column name={Basis}, 
        string type,
        column type=|c| 
    },       
    columns/Nx/.style={
        column name={$N_x$}, 
        int detect, 
        column type=c|
    },
    columns/ErrDet/.style={
        column name={$L^1$ Error}, 
        string replace={NaN}{},
        sci, sci zerofill, precision=2,
        column type=c|
    },  
    columns/OrdDet/.style={
        column name={Order}, 
        string replace={NaN}{},
        fixed, fixed zerofill, precision=2,
        column type=c|,
    }
]{graphics/data/NDG_check.dat}%
  \tikzexternaldisable%

      \caption{Convergence of the Nodal Discontinuous Galerkin (NDG) scheme: 
    The empirical order of accuracy closely aligns with the theoretical expectations for varying polynomial bases (NDG1, NDG2, and NDG3) as the spatial grid resolution ($N_x$) increases. 
    This optimal convergence validates the accuracy of the baseline solver implementation prior to introducing low-rank tensor approximations.}
    \label{fig:NDG_check}
\end{figure}

Next, we investigate the computational speed up obtained from the randomized sketching methods. We now assume a rank-6 initial condition
\begin{equation}
    f_0(x;\bxi)=\sum_{\alpha=1}^{6} g_{\alpha}(x;\bxi),
\end{equation}
where 
\begin{equation}
    g_{\alpha}(x;\bxi)=\frac{1+A_{\alpha}\sin(\omega_{\alpha}x)}{(2\pi)^{3/2}}\exp\left(-\frac{|\bxi-\bbeta_{\alpha}|^2}{2}\right).
\end{equation}
We took the values $(A_{\alpha},\omega_{\alpha},\bbeta_{\alpha})$, $1\leq\alpha\leq 6$ in order from the following sets:
\begin{align*}
A &= \left\{ \frac{1}{2}, \frac{1}{4}, \frac{3}{20}, \frac{1}{5}, \frac{9}{10}, \frac{1}{100} \right\}, \\
\omega &= \left\{ 1, \frac{1}{100}, 2, \frac{1}{5}, \frac{7}{10}, 10 \right\}, \\
\boldsymbol{\beta} &= \left\{ \begin{bmatrix} -1/20 \\ 0 \\ 0 \end{bmatrix}, \begin{bmatrix} 1/20 \\ 0 \\ 0 \end{bmatrix}, \begin{bmatrix} 0 \\ -3/40 \\ 0 \end{bmatrix}, \begin{bmatrix} 0 \\ 3/40 \\ 0 \end{bmatrix}, \begin{bmatrix} 0 \\ 0 \\ -1/40 \end{bmatrix}, \begin{bmatrix} 0 \\ 0 \\ 1/40 \end{bmatrix} \right\}.
\end{align*}

\Cref{fig:NDG_increasing_order} shows the computational speedup of \texttt{KronSum-Tucker} and \texttt{KRPSum-Tucker} relative to the (naive) \texttt{TuckerAxby}+\texttt{TuckerRounding} when increasing basis order. We set the spatial mesh $N_x=32$, parameter space mesh $N_{\bxi}=100$ with $\Xi=12$, tolerance $10^{-6}$, final time $T_f=0.005$, oversampling parameter $p=5$, and time-step size $\Delta t=\vartheta\Delta x$, where $\vartheta=0.1/((2k+1)\Xi)$. Letting the NDG order increase for $k=0,1,2,3$, we see that \texttt{KronSum-Tucker} and \texttt{KRPSum-Tucker} saw up to 5x speedup. Similarly, both randomized sketching methods observed great speedup when refining the velocity mesh $N_{\bxi}$. The results shown in \Cref{fig:NDG_increasing_Nv} use spatial mesh $N_x=32$, tolerance $10^{-6}$, final time $T_f=0.005$, oversampling parameter $p=5$, and NDG2 with $\Delta t=\vartheta\Delta x$. In both situations increasing either the size of the basis or the velocity mesh refinement, \texttt{KronSum-Tucker} and \texttt{KRPSum-Tucker} had almost comparable speedups. This is because unlike the \emph{cookie problem} where we had skewed modes, here the modes are not skewed since $N_{\xi_x}=N_{\xi_y}=N_{\xi_z}$; see \Cref{sec:Subrank_Selection} and \cite{MinLB24} (Subsection 4.3). Our results for this rank-6 initial condition suggest that for certain transport equations with low-rank solutions, the sketching-based techniques proposed in this paper become advantageous with higher-order bases and finer grids, both of which are often desired.

\begin{figure}[t]
    \centering    
    \begin{subfigure}{0.49\textwidth}
        \centering
  \tikzexternalenable%
  \tikzsetnextfilename{NDG_increase_order}%
  \begin{tikzpicture}[font = \plotfontsize\normalfont]
  \pgfplotstableread{graphics/data/NDG_increase_order.dat}\tableBasis

  \begin{axis}[
    width              = \linewidth,
    height             = 5cm,
    enlarge x limits   = 0.15,
    xmin               = 1,
    xmax               = 4,
    ymin               = 0,
    xtick              = {1, 2, 3, 4},
    ytick              = {1, 2, 3, 4, 5, 6},
    xlabel             = {NDG Order},
    ylabel             = {Speedup},
    ylabel style       = {yshift = -.3em},
    ymajorgrids        = true,
    scaled x ticks     = false,
  ]
    
    \addplot[Kron_Style_line,mark=square*, mark options={black}] 
    table[x=Basis, y=Spd_Kron] {\tableBasis};

    \addplot[KRP_Style_line,mark=*, mark options={black}] 
    table[x=Basis, y=Spd_KRP] {\tableBasis};

    \addplot[Naive_Style_line, mark=triangle*, mark options={black}] coordinates {(1,1) (2,1) (3,1) (4,1)};

  \end{axis}
\end{tikzpicture}%
  \tikzexternaldisable%

        \caption{Speed up for increasing order}
        \label{fig:NDG_increasing_order}
    \end{subfigure}
    \hfill
    \begin{subfigure}{0.49\textwidth}
        \centering
  \tikzexternalenable%
  \tikzsetnextfilename{NDG_increase_Nv}%
  \begin{tikzpicture}[font = \plotfontsize\normalfont]
  \pgfplotstableread{graphics/data/NDG_increase_Nv.dat}\tableNv

  \begin{axis}[
    width              = \linewidth,
    height             = 5cm,
    enlarge x limits   = 0.15,
    xmin               = 16,
    xmax               = 256,
    ymin               = 0,
    ymax               = 35,
    xmode              = log,
    log ticks with fixed point, 
    xtick              = {16, 32, 64, 128, 256},
    ytick              = {5, 10, 15, 20, 25, 30},
    xlabel             = {Number of Velocity Cells},
    ylabel             = {Speedup},
    ylabel style       = {yshift = -.3em},
    ymajorgrids        = true,
    scaled x ticks     = false,
  ]
    
    \addplot[Kron_Style_line,mark=square*, mark options={black}] 
    table[x=Nv, y=Spd_Kron] {\tableNv};

    \addplot[KRP_Style_line,mark=*, mark options={black} ]
    table[x=Nv, y=Spd_KRP] {\tableNv};

    \addplot[Naive_Style_line, mark=triangle*, mark options={black}] coordinates {(16,1) (32,1) (64,1) (128,1) (256,1)};

  \end{axis}
\end{tikzpicture}%
  \tikzexternaldisable%

        \caption{Speed up for Increasing $N_{\bxi}$}
        \label{fig:NDG_increasing_Nv}
    \end{subfigure}
    \vspace{.5\baselineskip}
  \tikzexternalenable%
  \tikzsetnextfilename{Legend}%
  \begin{tikzpicture}[font = \plotfontsize\normalfont]
  \begin{axis}[%
    hide axis, 
    scale only axis,
    width          = .8\linewidth,
    height         = .2\linewidth,
    xmin           = 0,
    xmax           = 1,
    ymin           = 0,
    ymax           = 1,
    legend columns = 3,    
    legend style   = {
      at     = {(0.5,0.5)}, 
      anchor = center,
      draw   = none,       
      /tikz/every even column/.append style = {column sep = 0.5cm}
    }
  ]
    
    
    \addlegendimage{Kron_Style_line, mark=square*, mark options={black}}
    \addlegendentry{\texttt{KronSum-Tucker}}

    \addlegendimage{KRP_Style_line,mark=*, mark options={black}}
    \addlegendentry{\texttt{KRPSum-Tucker}}

    \addlegendimage{Naive_Style_line,mark=triangle, mark options={black}}
    \addlegendentry{\texttt{TuckerAxby}+\texttt{TuckerRounding}}

  \end{axis}
\end{tikzpicture}%
  \tikzexternaldisable%

    \caption{Computational speedup of the randomized tensor approximations for the Nodal Discontinuous Galerkin (NDG) scheme. The efficiency of \texttt{KronSum-Tucker} and \texttt{KRPSum-Tucker} improves significantly as problem complexity grows relative to the \texttt{TuckerAxby+TuckerRounding}. (a) Increasing the NDG order yields up to a 5x speedup. (b) Refining the parameter space resolution demonstrates exponential performance gains, achieving nearly a 30x speedup at 256 velocity cells.}
    \label{fig:NDG_increasing_order_Nv}
  
\end{figure}

\section{Conclusions}
\label{sec:conclusions}

In this work, we introduced an efficient, sketching-based methodology for the summation of tensors in the Tucker format. We established the algebraic foundations of this approach, developing a framework that exploits the structure of Khatri-Rao and Kronecker products to perform arithmetic directly on compressed representations. This strategy successfully circumvents the explicit formation of massive intermediate core tensors. To ensure accurate low-rank approximations without the severe memory bottlenecks of continuous pairwise rounding, we also incorporated a spectral-based effective rank estimation method alongside a subrank selection heuristic. Finally, we supported our methodology with a comprehensive computational complexity analysis of the proposed algorithms.
The numerical examples verified that the proposed randomized approaches fundamentally resolve the intermediate rank swelling and dynamic range vulnerabilities associated with standard sequential truncation. Furthermore, evaluations on the parametric \emph{cookie problem} via \texttt{Tucker-GMRES} and a linear conservation law via a Nodal Discontinuous Galerkin scheme demonstrated that our sketching methods yield substantial computational speedups over current deterministic baselines while maintaining high accuracy. In contrast, omitting a rounding procedure entirely allows intermediate tensor ranks to grow unchecked after successive operations, leading to prohibitive memory consumption and rapidly escalating computational times that quickly render the simulations intractable.

While we motivated the use of the KRP and Kronecker sketching operators for structure preservation, the heuristic subrank selection presently relies on the exact construction of an energy-weighted Gram matrix, which scales with the number of summands. It is not yet proven whether an approximate or iterative sketching technique could estimate this effective rank more efficiently without compromising the error bounds of the final core projection.
Theoretical and algorithmic advancements on this topic will have strong implications for the application of randomized tensor arithmetic in other computationally demanding research areas, such as the real-time simulation of high-dimensional kinetic equations, nonlinear hyperbolic conservation laws, and dynamic model order reduction.

\section*{Acknowledgments}
MP acknowledges support from the National Science Foundation (NSF) under Grant No. DMS-2410699. Any opinions, findings, conclusions, or recommendations expressed in this material are those of the authors and do not necessarily reflect the views of the National Science Foundation. HAD acknowledges support from the Ada Lovelace Centre Programme at the Scientific Computing Department, STFC. GB acknowledges support from the US Department of Energy, Office of Science, Advanced Scientific Computing Research program under awards DE-SC-0023296 and DE-SC0025394. AGO gratefully acknowledges the support from the Oden Institute of Computational Engineering and Sciences. JQ acknowledges support from DoD MURI FA9550-24-1-0254 and the U.S. Department of Energy under grant DE-SC0023164. RS and MP acknowledge Advanced Research Computing at Virginia Tech for providing computational resources and technical support that have contributed to the results reported within this paper. URL: https://arc.vt.edu/. WT was supported by the DOE Office of Applied Scientific Computing Research (ASCR) through the Mathematical Multifaceted Integrated Capability Centers program. MP, AGO, JQ, JN, and WT acknowledge support from the NSF under Grant No. DMS-1929284 for a research stay at the Institute for Computational and Experimental Research in Mathematics in Providence, RI, during the \emph{Empowering a Diverse Computational Mathematics Research Community} program where preliminary discussions of this work took place.

\addcontentsline{toc}{section}{References}
\bibliographystyle{plainurl}
\bibliography{bibtex/myref}

\end{document}